\input amstex\documentstyle{amsppt}  
\pagewidth{12.5cm}\pageheight{19cm}\magnification\magstep1
\topmatter
\title{Unipotent classes and special Weyl group representations}\endtitle   
\author G. Lusztig\endauthor
\address{Department of Mathematics, M.I.T., Cambridge, MA 02139}\endaddress  
\thanks{Supported in part by the National Science Foundation}\endthanks
\endtopmatter   
\document

\define\du{\dot u}

\define\uG{\un G}

\define\hx{\hat x}

\define\da{\dagger}

\define\si{\sim}

\define\sqc{\sqcup}

\define\tca{\ti\ca}

\define\tiz{\ti\z}

\define\lb{\linebreak}

\redefine\sp{\spadesuit}
\define\part{\partial}
\define\em{\emptyset}

\define\ra{\rangle}
\define\n{\notin}
\define\iy{\infty}
\define\m{\mapsto}
\define\do{\dots}
\define\la{\langle}

\define\lra{\leftrightarrow}

\define\sm{\smallmatrix}
\define\esm{\endsmallmatrix}
\define\sub{\subset}    
\define\bxt{\boxtimes}
\define\T{\times}
\define\ti{\tilde}
\define\nl{\newline}
\redefine\i{^{-1}}
\define\fra{\frac}
\define\un{\underline}

\define\Ad{\text{\rm Ad}}

\define\Ind{\text{\rm Ind}}

\define\a{\alpha}
\redefine\b{\beta}

\redefine\d{\delta}
\define\e{\epsilon}

\define\io{\iota}
\redefine\o{\omega}

\define\r{\rho}

\redefine\t{\tau}

\define\k{\kappa}
\redefine\l{\lambda}
\define\z{\zeta}
\define\x{\xi}

\define\Om{\Omega}

\define\bb{\bold b}

\define\qq{\bold q}

\define\zz{\bold z}

\define\CC{\bold C}

\define\NN{\bold N}

\define\QQ{\bold Q}

\define\WW{\bold W}
\define\ZZ{\bold Z}

\define\ca{\Cal A}

\define\ce{\Cal E}

\define\cg{\Cal G}

\define\ci{\Cal I}
\define\cj{\Cal J}

\define\cp{\Cal P}

\define\car{\Cal R}
\define\cs{\Cal S}

\define\cz{\Cal Z}
\define\cx{\Cal X}

\define\fa{\frak a}

\define\fc{\frak c}

\define\fI{\frak I}

\define\fS{\frak S}

\define\te{\ti e}

\define\tp{\ti p}

\define\tz{\ti z}
\define\tx{\ti x}
\define\ty{\ti y}

\define\tE{\ti E}

\define\tI{\ti I}
\define\tJ{\ti J}

\define\tS{\ti S}

\define\tX{\ti X}
\define\tY{\ti Y}

\define\tiO{\ti\Om}

\define\Irr{\text{\rm Irr}}

\define\che{\check}
\define\cha{\che{\a}}

\define\tib{\ti\b}

\define\uzz{\un\zz}
\define\tzz{\ti\zz}
\define\bcs{\bar{\cs}}
\define\tcs{\ti\cs}
\define\ddu{\ddot u}
\define\AL{Al}
\define\ALL{AL}
\define\SPE{L1}
\define\CL{L2}
\define\UC{L3}
\define\OR{L4}
\define\ICC{L5}
\define\USU{L6}
\define\SH{S1}
\define\SHO{S2}
\define\SP{Sp}

\head Introduction \endhead
\subhead 0.1\endsubhead
Let $G$ be a simple adjoint algebraic group over $\CC$ and let $\cx$ be the set of unipotent conjugacy classes in
$G$. Let $C\in\cx$ and let $u\in C$. The following invariants of $C$ are important in representation theory:

-the dimension $\bb_C$ of the fixed point set of $\Ad(u)$ on the flag manifold of $G$;

-the number $\zz_C$ of connected components of the centralizer of $u$ in $G$;

-the number $\tzz_C$ of connected components of the centralizer of a unipotent element in the simply connected
covering of $G$ which projects to $u$;

-the irreducible representation $\r_C$ of the Weyl group $\WW$ of $G$ corresponding to $C$ and the constant local
system under the Springer correspondence \cite{\SP}.
\nl
Let $\tcs_\WW$ be the set of isomorphism classes of irreducible representations of $\WW$ of the form $\r_C$ for
some $C\in\cx$. It is known \cite{\SP} 
that $C\m\r_C$ is a bijection $\cx@>\si>>\tcs_\WW$.

Note that the definition of each of $\bb_C$, $\zz_C$, $\tzz_C$ is based on considerations of algebraic geometry 
and in the case of $\tcs_\WW$, also on considerations of \'etale cohomology.

In \cite{\SPE, Sec.9} I conjectured that $\tcs_\WW$, $C\m\bb_C$ and $C\m\zz_C$ can be determined purely in terms 
of data involving the Weyl group $\WW$ (more precisely, the "special representations" of the "parahoric" subgroups
of $\WW$, see 1.1, 1.2). At that time I could only prove this conjecture for $\tcs_\WW$ and for $C\m\bb_C$ assuming that $G$ is
of classical type (my proof was based on \cite{\SH}) and a little later for $G$ of type $F_4$ (based on 
\cite{\SHO}). In \cite{\ALL} the conjecture for $\tcs_\WW$ and $C\m\bb_C$ was established for $G$ of type 
$E_6,E_7,E_8$. At the time \cite{\OR} was written, I proved the remaining 
conjecture of \cite{\SPE} (concerning $C\m\zz_C$); this was stated in 
\cite{\OR, 13.3}. For
classical groups the proof involved a new description (in terms of "symbols") of the Springer correspondence for 
classical groups (given in \cite{\ICC}) while for exceptional groups this was a purely mechanical verification 
based on the tables \cite{\AL}. The conjecture of \cite{\SPE} is restated and proved here as Theorem 
1.5(a),(b1),(b2). At the same time we state and prove a complement to that 
conjecture, namely that $C\m\tzz_C$ is 
determined purely in terms of data involving $\WW$ (see Theorem 
1.5(b3)). Note that for classical groups this involves some combinatorial considerations while for exceptional 
groups this involves only a purely mechanical verification based on the known tables.

{\it Notation.} For a finite set $F$ let $|F|$ be the cardinal of $F$. For $i,j$ in $\ZZ$ we set 
$[i,j]=\{n\in\ZZ;i\le n\le j\}$. For $x,y$ in $\ZZ$ we write $x\ll y$ if $x\le y-2$.

\head Contents\endhead
1. Statement of the main result.

2. Combinatorics.

3. Type $A_{n-1}$.

4. Type $B_n$.

5. Type $C_n$.

6. Type $D_n$.

7. Exceptional types.

Index.

\head 1. Statement of the main result\endhead
\subhead 1.1\endsubhead
Let $W$ be a finite crystallographic Coxeter group. Let $\Irr(W)$ be the set of isomorphism classes of irreducible
representations of $W$ over $\QQ$. If $E\in\Irr(W)$ and $E'$ is a finite dimensional $\QQ[W]$-module, let 
$[E:E']_W$ be the multiplicity of $E$ in $E'$. Let $S^i_W$ be the $i$-th symmetric power of the reflection 
representation of $W$. For any $E\in\Irr(W)$ we define integers $f_E\ge1$, $a_E\ge0$ by the requirement that the 
generic degree of the Hecke algebra representation corresponding to $E$ is of the form 
$\fra{1}{f_E}\qq^{a_E}+$higher powers of $\qq$ ($\qq$ is an indeterminate); let $b_E$ be the smallest integer 
$i\ge0$ such that $[E:S^i_W]_W\ge1$. As observed in \cite{\SPE, Sec.2}, we have $a_E\le b_E$ for any $E\in\Irr(W)$; 
following \cite{\SPE, Sec.2} we set $\cs_W=\{E\in\Irr(W);a_E=b_E\}$; this is the set of "special representations" of
$W$. Let $\Irr(W)^\da=\{E\in\Irr(W);[E:S^{b_E}_W]_W=1\}$. We have $\cs_W\sub\Irr(W)^\da$. 

\subhead 1.2\endsubhead
In this paper we fix a root datum of finite type $\car=(Y,X,\cha_i,\a_i(i\in I),\la,\ra)$. (Here $Y,X$ are free 
abelian groups of finite rank, $\la,\ra:Y\T X@>>>\ZZ$ is a perfect pairing, $\cha_i\in Y$ are the simple coroots 
and $\a_i\in X$ are the simple roots.) We assume that $I\ne\em$ and that $\car$ is of adjoint type that is, 
$\{\a_i;i\in I\}$ is a $\ZZ$-basis of $X$. Let $R\sub X$ (resp. $\che R\sub Y$) be the set of roots (resp. 
coroots); let $\cha\lra\a$ be the canonical bijection $\che R\lra R$. We assume that $\car$ is irreducible that 
is, there is a unique $\a_0\in R$ such that $\cha_0-\cha_i\n\che R$ for any $i\in I$. Let $\tI=I\sqc\{0\}$. For 
$i\in\tI$ let $s_i:X@>>>X$ be the reflection determined by $\a_i,\cha_i$. Let $\WW$ be the subgroup of $GL(X)$ 
generated by $\{s_i;i\in I\}$, a finite crystallographic Coxeter group containing $s_0$. The elements 
$s_i(i\in\tI)$ in $\WW$ satisfy the relations of the affine Weyl group of type dual to that of $\car$. Let 
$\tca=\{J;J\subsetneqq\tI\}$. For any $J\in\tca$, let $\WW_J$ be the subgroup of $\WW$ generated by 
$\{s_i;i\in J\}$, a finite crystallographic Coxeter group with set of generators $\{s_i;i\in J\}$, said to be a 
{\it parahoric subgroup} of $\WW$. 

Let $\Om$ be the (commutative) subgroup of $\WW$ consisting of all $\o\in\WW$ 
such that $\o(\a_i)=\a_{\un\o(i)}$ ($i\in\tI$) for some (necessarily unique) 
permutation $\un\o:\tI@>\si>>\tI$.

\subhead 1.3\endsubhead
If $J\in\tca$ and $E_1\in\Irr(\WW_J)^\da$, there is a unique $E\in\Irr(\WW)$ such that $b_E=b_{E_1}$ and 
$[E:\Ind_{\WW_J}^\WW(E_1)]_\WW\ge1$. (Then $[E:\Ind_{\WW_J}^\WW(E_1)]_\WW=1$ and $E\in\Irr(\WW)^\da$.) We write 
$E=j_{\WW_J}^\WW(E_1)$. Let $E\in\Irr(\WW)$ and let
$$\cz_E=\{(J,E_1);J\in\tca,E_1\in\cs_{\WW_J},E=j_{\WW_J}^\WW(E_1)\}.$$
Let 
$$\bcs_\WW=\{E\in\Irr(\WW);\cz_E\ne\em\}.$$
Let $E\in\bcs_\WW$. We set 
$$\fa_E=\max_{(J,E_1)\in\cz_E}f_{E_1}.$$
Let $\cz_E^\sp=\{(J,E_1)\in\cz_E;f_{E_1}=\fa_E\}$. We have $\cz_E^\sp\ne\em$. 

If $(J,E_1)\in\cz_E$ and $\o\in\Om$ then $\Ad(\o):\WW_J@>\si>>\WW_{\un\o(J)}$ carries $E_1$ to a representation 
${}^\o E_1\in\cs_{\WW_{\un\o(J)}}$ such that $\Ind_{\WW_J^\WW}(E_1)=\Ind_{\WW_{\un\o(J)}^\WW}({}^\o E_1)$, 
$b_{{}^\o E_1}=b_{E_1}$ and $f_{{}^\o E_1}=f_{E_1}$. It follows that 
$j_{\WW_J}^\WW(E_1)=j_{\WW_{\un\o(J)}}^\WW({}^\o E_1)$. Thus $(\un\o(J),{}^\o E_1)\in\cz_E$ and
$\o:(J,E_1)\m(\un\o(J),{}^\o E_1)$ is an action of $\Om$ on $\cz_E$. This restricts to an action of $\Om$ on 
$\cz_E^\sp$. The stabilizer in $\Om$ of $(J,E_1)\in\cz_E^\sp$ for this action is denoted by $\Om_{J,E_1}$. We set 
$$\fc_E=\max_{(J,E_1)\in\cz_E^\sp}|\Om_{J,E_1}|.$$

\subhead 1.4\endsubhead
Let $G$ be a semisimple (adjoint) algebraic group over $\CC$ with root datum $\car$. 
Let $\cx$, $C\m\r_C$, $C\m\bb_C$, $C\m\zz_C$, $C\m\tzz_C$, $\tcs_\WW$ be as in 0.1.

\proclaim{Theorem 1.5}(a) $\tcs_\WW=\bcs_\WW$. 

(b) Let $C\in\cx$. Set $E=\r_C\in\bcs_\WW$. Then:

(b1) $\bb_C=b_E$;

(b2) $\zz_C=\fa_E$;

(b3) $\tzz_C/\zz_C=\fc_E$.
\endproclaim
For exceptional types the proof of (a),(b1)-(b3) consists in examining the existing tables. Some relevant data is
collected in \S7. The proof for the classical types is given in \S3-\S6 after combinatorial preliminaries in
1.9-1.11 and \S2.

\subhead 1.6\endsubhead
Let $G'$ be a connected reductive group over $\CC$  such that $G$ is the quotient of $G'$ by its centre. 

Note that 1.5(a) is closely connected to the definition of a unipotent support
of a character sheaf on $G'$ provided by \cite{\USU, 10.7}. In fact, 
\cite{\USU, 10.7(iii)} provides a proof of the inclusion $\bcs_\WW\sub\tcs_\WW$
without case by case checking.

For any $g\in G'$ let $g_u$ be the unipotent part of $g$. We now state an 
alternative conjectural definition of the unipotent support of a character
sheaf on $G'$.

\proclaim{Conjecture 1.7} Let $A$ be a character sheaf on $G'$. There exists a unique unipotent class $C$ in $G'$
such that:

(i) $A|_{\{g\}}\ne0$ for some $g\in G'$ with $g_u\in C$;

(ii) if $g'\in G'$ satisfies $A|_{\{g'\}}\ne0$ then the conjugacy class of $g'_u$ in $G'$ has dimension $<\dim(C)$.
\endproclaim

\subhead 1.8\endsubhead
Theorem 1.5 remains valid if $\CC$ is replaced by an algebraically closed field whose characteristic is either $0$
or a prime which is good for $G$ and which (if $G$ is of type $A_{n-1}$) does not divide $n$.

\subhead 1.9\endsubhead
In the rest of this section we discuss some preliminaries to the proof of 1.5.

If $J,J'\in\tca$, $J\sub J'$ and $E_1\in\Irr(\WW_J)^\da$, there is a unique $E'_1\in\Irr(\WW_{J'})$ such that 
$b_{E_1}=b_{E'_1}$ and $[E'_1:\Ind_{\WW_J}^{\WW_{J'}}(E_1)]_{\WW_{J'}}\ge1$. (Then 
$[E'_1:\Ind_{\WW_J}^{\WW_{J'}}(E_1)]_{\WW_{J'}}=1$ and $E'_1\in\Irr(\WW_{J'})^\da$.) We write 
$E'_1=j_{\WW_J}^{\WW_{J'}}(E_1)$. Note that 

(a) $j_{\WW_J}^\WW(E_1)=j_{\WW_{J'}}^\WW(j_{\WW_J}^{\WW_{J'}}(E_1))$;

(b) if, in addition, $E_1\in\cs_{\WW_J}$, then $E'_1\in\cs_{\WW_{J'}}$ and $f_{E_1}\le f_{E'_1}$.
\nl
(See \cite{\SPE, Sec.4}.)

Let $\cp'$ be the collection of parahoric subgroups $W$ of $\WW$ such that $W=\WW_J$ for some $J\sub\tI$, 
$|J|=|\tI|-1$. From (a),(b) we see that
$$\bcs_\WW=\{E\in\Irr(\WW);E=j_W^\WW(E_1)\text{ for some }W\in\cp'
\text{ and some }E_1\in\cs_W\},\tag c$$
$$\fa_E=\max_{(J,E_1)\in\cz_E;|J|=|\tI|-1}f_{E_1}\text{ for }E\in\bcs_\WW.\tag d$$
If $W=W_1\T W_2$ where $W_1,W_2$ are finite crystallographic Coxeter groups and $E_1\in\Irr(W_1)$, 
$E_2\in\Irr(W_2)$ then $E:=E_1\bxt E_2\in\Irr(W)$ belongs to $\cs_W$ if and only if $E_1\in\cs_{W_1}$ and 
$E_2\in\cs_{W_2}$; in this case we have 

(e) $a_E=a_{E_1}+a_{E_2}$, $f_E=f_{E_1}f_{E_2}$. 

\subhead 1.10\endsubhead
We show: 

(a) {\it if $J,J'\in\tca$ and $\WW_J=\WW_{J'}\ne\WW$ then $J=J'$.}
\nl
It is enough to show that if $J,J'\in\tca$ and $\WW_J\sub\WW_{J'}\ne\WW$ then $J\sub J'$. To see this we may
assume that $J$ consists of a single element $j$. We have $s_j\in\WW_{J'}$. Assume that $j\n J'$. If 
$J'\cup\{j\}\ne\tI$ then $\WW_{J'\cup\{j\}}$ is a Coxeter group on the generators $\{s_h;h\in J'\cup\{j\}\}$. In 
particular $s_j$ is not contained in the subgroup $\WW_{J'}$ generated by $\{s_h;h\in J'\}$, a contradiction. Thus
we have $J'\cup\{j\}=\tJ$. We see that $\WW_{J'}$ contains $\{s_h;h\in J'\cup\{j\}\}$ which generates $\WW$. Thus
$\WW_J=\WW$ which is again a contradiction. This proves (a).

\subhead 1.11\endsubhead
For a subgroup $\tiO$ of $\Om$ let $\cp^{\tiO}$ be the collection of parahoric subgroups $W$ of $\WW$ such that 
$W=\WW_J$ for some $J\in\tca$ where $J$ is $\tiO$-stable and is maximal with this property. From the
definitions we have
$$\fc_E=\max|\tiO|,$$
where the maximum is taken over all subgroups
 $\tiO\sub\Om$ and all $(J,E_1)\in\cz_E^\sp$ such that $\WW_J\in\cp^{\tiO}$,
$\tiO\sub\Om_{J,E_1}$.

\head 2. Combinatorics\endhead
\subhead 2.1\endsubhead
In this section we fix $m\in\NN$.

Let $Z_m=\{z_*=(z_0,z_1,z_2,\do,z_m)\in\NN^{m+1};z_0<z_1<\do<z_m\}$. Let $z_*^0=z_*^{0,m}=(0,1,2,\do,m)\in Z_m$. 
For any $z_*\in Z_m$ we have $z_*-z^0_*\in\NN^{m+1}$. Hence 

$\r_0:Z_m@>>>\NN$, $z_*\m\sum_{i\in[0,m]}(z_i-z^0_i)$ and 

$\b_0:Z_m@>>>\NN$, $z_*\m\sum_{0\le i<j\le m}(z_i-z_i^0)$
\nl
are well defined. For any $n\in\NN$ we set $Z_m^n=\{z_*\in Z_m;\r_0(z_*)=n\}$. 

\subhead 2.2\endsubhead
Let $X_m$ be the set of all $x_*=(x_0,x_1,x_2,\do,x_m)\in\NN^{m+1}$ such that $x_i\le x_{i+1}$ for $i\in[0,m-1]$,
$x_i<x_{i+2}$ for $i\in[0,m-2]$. For $x_*\in X_m$ let $\fS(x_*)$ be the set of all $i\in[0,m]$ such that 
$x_{i-1}<x_i<x_{i+1}$ (with the convention $x_{-1}=-\iy,x_{m+1}=\iy$). Note that 

(a) $|\fS(x_*)|\cong m-1\mod2$;

(b) $\fS(x_*)=\em$ if and only if $m$ is odd and $x_i=x_{i+1}$ for $i=0,1,\do,(m-1)/2$. 

\subhead 2.3\endsubhead
Let $Y_m$ be the set of all $y_*=(y_0,y_1,y_2,\do,y_m)\in\NN^{m+1}$ such that $y_i\le y_{i+1}$ for $i\in[0,m-1]$,
$y_i\ll y_{i+2}$ for $i\in[0,m-2]$. For $y_*\in Y_m$ let $\fI(y_*)$ be the set of all intervals $[i,j]\sub[0,m]$
(with $i\le j$) such that 

$y_{i-1}-(i-1)<y_i-i=y_{i+1}-{i+1}=\do=y_j-j<y_{j+1}-(j+1)$
\nl
(with the convention $y_{-1}=-\iy,y_{m+1}=\iy$). We have 

(a) $\fI(y_*)=\em$ if and only if $m$ is odd and $y_i=y_{i+1}$ for $i=0,1,\do,(m-1)/2$.
\nl
Let 

$\fI'(y_*)=\{\ci\in\fI(y_*);\ci=[i,j]\text{ with $|[i,j]|=$odd}\}$, 

$\fI''(y_*)=\{\ci\in\fI(y_*);\ci=[i,j]\text{ with $|[i,j]|=$even}\}$. 
\nl
We have 

(b) $|\fI'(y_*)|\cong m-1\mod2$.
\nl
Let $R(y_*)$ be the set of all $k\in[0,m]$ such that $k=i$ or $k=j$ for some (necessarily unique) 
$[i,j]\in\fI(y_*)$. Let $R_0(y_*)$ be the set of all $k\in[0,m]$ such that $k=i$ for some (necessarily unique) 
$[i,j]\in\fI(y_*)$ with $i=j$. Clearly,

(c) $|R(y_*)|+|R_0(y_*)|=2|\fI(y_*)|.$

\subhead 2.4\endsubhead
Let $x_*,x'_*\in X_m$ and let $y_*=x_*+x'_*\in\NN^{m+1}$. Note that $y\in Y_m$. If $k\in\fS(x_*)$ then 
$x_{k-1}<x_k<x_{k+1}$, $x'_{k-1}\le x'_k\le x'_{k+1}$ (and at least one of the last two $\le$ is $<$). Hence 
$y_{k-1}<y_k<y_{k+1}$ (and at least one of the last two $<$ is $\ll$). Hence $k\in R(y_*)$. Thus 
$\fS(x_*)\sub R(y_*)$. Similarly, $\fS(x'_*)\sub R(y_*)$. We see that $\fS(x_*)\cup\fS(x'_*)\sub R(y_*)$. If 
$k\in\fS(x_*)\cap\fS(x'_*)$ then $x_{k-1}<x_k<x_{k+1}$, $x'_{k-1}<x'_k<x'_{k+1}$ hence $y_{k-1}\ll y_k\ll y_{k+1}$
so that $k\in R_0(y_*)$. Thus, $\fS(x_*)\cap\fS(x'_*)\sub R_0(y_*)$ and
$$|\fS(x_*)|+|\fS(x'_*)|=|\fS(x_*)\cup\fS(x'_*)|+|\fS(x_*)\cap\fS(x'_*)|\le|R(y_*)|+|R_0(y_*)|.$$
Using this and 2.3(c) we see that
$$|\fS(x_*)|+|\fS(x'_*)|\le2|\fI(y_*)|,\tag a$$
with equality if and only if $\fS(x_*)\cup\fS(x'_*)=R(y_*)$ and $\fS(x_*)\cap\fS(x'_*)=R_0(y_*)$.
 
\subhead 2.5\endsubhead
Let $y_*\in Y_m$. We consider a partition $[0,m]=\cj_0\sqc\cj_1\sqc\do\cj_t$ where for each $s\in[0,t]$ we have
$\cj_s=[m_s,m'_{s+1}]$ with $m_s\le m'_{s+1}$, $m_0=0$, $m'_{t+1}=m$ and for each $s\in[1,t]$ we have 
$m_s=m'_s+1$. We require that for $s\in[1,t]$ we have $y_{m'_s}\ll y_{m_s}$ and for any $s\in[0,t]$ we have either

(i) $|\cj_s|=2$ and $(y_{m_s},y_{m'_{s+1}})=(a_s,a_s)$, or

(ii) $(y_{m_s},y_{m_s+1},\do,y_{m'_{s+1}})=(a_s,a_s+1,a_s+2,\do)$. 
\nl
for some $a_s\in\NN$. Such a partition exists and is unique. Let

$\cg_1(y_*)=\{s\in[0,t];s\text{ is as in (i)}\}$, $\cg_2(y_*)=\{s\in[0,t];s\text{ is as in (ii)}\}$. 
\nl
We have 

$\fI(y_*)=\{[i,j];i=m_s,j=m'_{s+1}\text{ for some $s\in\cg_2(y_*)$}\}$;

$R(y_*)=\{i\in[0,m];i=m_s\text{ or $i=m'_{s+1}$ for some $s\in\cg_2(y_*)$}\}$;

$R_0(y_*)=\{i\in[0,m];i=m_s=m'_{s+1}\text{ for some $s\in\cg_2(y_*)$}\}$.

\subhead 2.6\endsubhead
Let $y_*\in Y_m$. Let $S'(y_*)$ be the set consisting of all pairs 

$x_*=(x_0,x_1,\do,x_m)$, $x'_*=(x'_0,x'_1,\do,x'_m)$
\nl
in $\NN^{m+1}$ which satisfy (i)-(iv) below (notation in 2.5):

(i) for any $s\in\cg_1(y_*)$ we have $(x_{m_s},x_{m'_{s+1}})=(u_s,u_s)$, $(x'_{m_s},x'_{m'_{s+1}})=(u'_s,u'_s)$, 
$u_s+u'_s=a_s$;
 
(ii) for any $s\in\cg_2(y_*)$ we have either

(ii1) $(x_{m_s},x_{m_s+1},\do,x_{m'_{s+1}})=(u_s,u_s+1,u_s+1,u_s+2,u_s+2,u_s+3,\do)$,
$(x'_{m_s},x'_{m_s+1},\do,x'_{m'_{s+1}})=(u'_s,u'_s,u'_s+1,u'_s+1,u'_s+2,u'_s+2,\do)$, $u_s+u'_s=a_s$, or

(ii2) $(x_{m_s},x_{m_s+1},\do,x_{m'_{s+1}})=(u_s,u_s,u_s+1,u_s+1,u_s+2,u_s+2,\do)$,

$(x'_{m_s},x'_{m_s+1},\do,x'_{m'_{s+1}})=(u'_s,u'_s+1,u'_s+1,u'_s+2,u'_s+2,u'_s+3,\do)$, $u_s+u'_s=a_s$;

(iii) for any $s\in[1,t]$ we have $x_{m'_s}<x_{m_s}$, $x'_{m'_s}<x'_{m_s}$;

(iv) if $\fI'(y_*)=\em$ then for any $s\in\cg_2(y_*)$, 
$(x_{m_s},x_{m_s+1},\do,x_{m'_{s+1}})$, \lb
$(x'_{m_s},x'_{m_s+1},\do,x'_{m'_{s+1}})$ are as in (ii1).
\nl
An element $(x_*,x'_*)$ of $S'(y_*)$ can be constructed by induction as follows. Assume that the entries 
$x_i,x'_i$ have been already chosen for $i\in\cj_0\cup\cj_1\cup\do\cup\cj_{s-1}$ for some $s\in[0,t]$ so that 
(i)-(iii) hold as far as it makes sense. In the case where $s>0$ let $\x=x_{m'_s},\x'=x'_{m'_s}$; in the case 
where $s=0$ let $\x=\x'=-\iy$. In any case we have $\x+\x'\le a_s-2$ hence we can find $u_s,u'_s$ in $\NN$ such 
that $\x<u_s$, $\x'<u'_s$, $u_s+u'_s=a_s$. (The number of choices is $y_{m_s}-y_{m'_s}-1$ if $s>0$ and $y_0+1$ if
$s=0$.) Then we define 

$(x_{m_s},x_{m_s+1},\do,x_{m'_{s+1}})$, $(x'_{m_s},x'_{m_s+1},\do,x'_{m'_{s+1}})$
\nl
by (i)
if $s\in\cg_1(y_*)$ and by (ii) if $s\in\cg_2(y_*)$. This gives two choices for each $s\in\cg_2(y_*)$ such that 
$|\cj_s|>1$, unless $\fI'(y_*)=\em$ when there is only one choice. This completes the inductive definition of 
$x_*,x'_*$. We see that $S'(y_*)\ne\em$.

Let $S(y_*)$ be the set of all $(x_*,x'_*)\in X_m\T X_m$ such that (v),(vi),(vii) below hold:

(v) $x_*+x'_*=y_*$,

(vi) $\fS(x_*)\cup\fS(x'_*)=R(y_*)$, $\fS(x_*)\cap\fS(x'_*)=R_0(y_*)$ (or equivalently 
$|\fS(x_*)|+|\fS(x'_*)|=2|\fI(y_*)|$),

(vii) if $\fI'(y_*)=\em$ (so that $m$ is odd), then $\fS(x'_*)=\em$.
\nl
From the definitions we see that $S(y_*)=S'(y_*)$. Hence 

(a) $S(y_*)\ne\em$.
\nl
From 2.4(a) we see that:

(b) if $\fI(y_*)=\em$ and $(x_*,x'_*)\in S(y_*)$ then $\fS(x_*)=\em$, $\fS(x'_*)=\em$.
\nl
On the other hand,

(c) if $\fI'(y_*)\ne\em$ and $(x_*,x'_*)\in S(y_*)$ then $\fS(x_*)\ne\em$, $\fS(x'_*)\ne\em$.
\nl
Indeed, let $[i,j]\in\fI'(y_*)$. Then we have either $i\in\fS(x_*),j\in\fS(x'_*)$ or $i\in\fS(x'_*),j\in\fS(x_*)$;
in both cases the conclusion of (c) holds.

\subhead 2.7\endsubhead
In this subsection we assume that $m$ is even, $\ge2$. We set 

$\tX_m=\{x_*\in X_m;x_0=0,x_1\ge1\}$, 
$\tY_m=\{y_*\in Y_m;y_1\ge1\}$. 
\nl
If $x_*\in X_m$, $x'_*\in\tX_m$, then $x_*+x'_*\in\tY_m$. 

Let $y_*\in\tY_m$ be such that 

(a) $y_0=0,y_1=1$. 
\nl
(Thus $\fI(y_*)$ contains an interval of form $[0,\a]$ hence $\fI(y_*)\ne\em$.) Let $\tS'(y_*)$ be the set 
consisting of all pairs $x_*=(x_0,x_1,\do,x_m)$, $x'_*=(x'_0,x'_1,\do,x'_m)$ in $\NN^{m+1}$ which satisfy the 
conditions (i)-(iii) in 2.6 together with conditions (i),(ii) below (notation in 2.5):

(i) for $s=0$ (necessarily in $\cg_2(y_*)$) we have 

$(x_0,x_1,\do,x_{m'_1})=(0,0,1,1,2,2,\do)$, $(x'_0,x'_1,\do,x'_{m'_1})=(0,1,1,2,2,3,3,\do)$
\nl
(so that $0\in\fS(x'_*)$);

(ii) if $\fI(y_*)=\{[0,\a]\}\cup\fI''(y_*)$ (so that $\fI'(y_*)=\{[0,\a]\}$) then 
for any $s\in\cg_2(y_*)-\{0\}$, $(x_{m_s},x_{m_s+1},\do,x_{m'_{s+1}})$, 
$(x'_{m_s},x'_{m_s+1},\do,x'_{m'_{s+1}})$ are as in 2.6(ii1).
\nl
We can construct an element in $\tS'(y_*)$ by the same method as in 2.6. In particular, $\tS'(y_*)\ne\em$. 

Now let $\tS(y_*)$ be the set of all $(x_*,x'_*)\in X_m\T\tX_m$ such that 

(iii) $x_*+x'_*=y_*$,

(iv) $\fS(x_*)\cup\fS(x'_*)=R(y_*)$, $\fS(x_*)\cap\fS(x'_*)=R_0(y_*)$ (or equivalently 
$|\fS(x_*)|+|\fS(x'_*)|=2|\fI(y_*)|$),

(v) if $\fI(y_*)=\{[0,\a]\}\cup\fI''(y_*)$, then $\fS(x'_*)=\{0\}$.
\nl
From the definitions we see that $\tS(y_*)=\tS'(y_*)$. Hence 

(b) $\tS(y)\ne\em$.
\nl
Note that

(c) if $\fI(y_*)=\{[0,\a]\}$ and $(x_*,x'_*)\in\tS(y_*)$, then $\fS(x_*)=\{\a\}$, $\fS(x'_*)=\{0\}$.
\nl
Indeed from 2.4(a) we see that $|\fS(x_*)|+|\fS(x'_*)|\le2$. On the other hand, we have $0\in\fS(x'_*)$ and
$\a\in\fS(x_*)$ (see (i)) since in this case $\a$ is even; (c) follows. Note that

(d) if $\fI'(y_*)$ contains at least one interval $\ne[0,\a]$ and $(x_*,x'_*)\in\tS(y_*)$, then $|\fS(x'_*)|\ge3$.
\nl
Indeed, let $[i,j]\in\fI'(y_*)$, $[i,j]\ne[0,\a]$. Then we have either $i\in\fS(x_*),j\in\fS(x'_*)$ or 
$i\in\fS(x'_*),j\in\fS(x_*)$. Since $0\in\fS(x'_*)$ we see that $|\fS(x'_*)|\ge2$. Since $|\fS(x'_*)|$ is odd we 
see that $|\fS(x'_*)|\ge3$.

\subhead 2.8\endsubhead
Let $x_*^0\in X_m$ be $(0,0,1,1,\do,(n-1),(n-1),n)$ if $m=2n$ and \lb $(0,0,1,1,\do,n,n)$ if $m=2n+1$. For any 
$x_*\in X_m$ we have $x_i\ge x^0_i$ for all $i\in[0,m]$. Hence 

$\r:X_m@>>>\NN$, $\x_*\m\sum_{i\in[0,m]}(x_i-x^0_i)$ and 

$\b:X_m@>>>\NN$, $x_*\m\sum_{0\le i<j\le m}(x_i-x_i^0)$
\nl
are well defined. 

Let $y_*^0\in Y_m$ be $(0,0,2,2,\do,(m-2),(m-2),m)$ if $m$ is even and \lb
$(0,0,2,2,\do,(m-1),(m-1))$ if $m$ is odd.
For any $y_*\in Y_m$ we have $y_i\ge y^0_i$ for all $i\in[0,m]$. Hence 

$\r':Y_m@>>>\NN$, $y_*\m\sum_{i\in[0,m]}(y_i-y^0_i)$ and 

$\b':Y_m@>>>\NN$, $y_*\m\sum_{0\le i<j\le m}(y_i-y_i^0)$
\nl
are
well defined. Since $x^0_*+x^0_*=y^0_*$ we have 

$\r'(x_*+x'_*)=\r(x_*)+\r(x'_*)$, $\b'(x_*+x'_*)=\b(x_*)+\b(x'_*)$
\nl
for any $x_*,x'_*\in X_m$. For any $n\in\NN$ we set $X_m^n=\{x_*\in X_m;\r(x_*)=n\}$, 
$Y_m^n=\{y_*\in Y_m;\r'(y_*)=n\}$.

Assume that $m=2k$, $k\ge1$. Let $\tx_*^0\in\tX_m$ be $(0,1,1,\do,k,k)$. For any $x_*\in\tX_m$ we have 
$x_i\ge\tx^0_i$ for all $i$. Hence 

$\ti\r:\tX_m@>>>\NN$, $\x_*\m\sum_{i\in[0,m]}(x_i-\tx^0_i)$ and

$\tib:\tX_m@>>>\NN$, $x_*\m\sum_{0\le i<j\le m}(x_i-\tx_i^0)$
\nl
are well defined. Let 
$\ty_*^0=(0,1,2,3,\do,m)\in Y_m$. For any $y_*\in\tY_m$ we have $y_i\ge\ty^0_i$ for all $i$. Hence 

$\ti\r':\tY_m@>>>\NN$, $y_*\m\sum_{i\in[0,m]}(y_i-\ty^0_i)$ and 

$\tib':\tY_m@>>>\NN$, $y_*\m\sum_{0\le i<j\le m}(y_i-\ty_i^0)$
\nl
are well defined. Since $x^0_*+\tx^0_*=\ty^0_*$ we have 

$\ti\r'(x_*+x'_*)=\r(x_*)+\ti\r(x'_*)$, $\tib'(x_*+x'_*)=\b(x_*)+\tib(x'_*)$ 
\nl
for any $x_*\in X_m$, $x'_*\in\tX_m$. For any $n\in\NN$ we set 

$\tX_m^n=\{x_*\in\tX_m;\ti\r(x_*)=n\}$, $\tY_m^n=\{y_*\in\tY_m;\ti\r'(y_*)=n\}$.

\subhead 2.9\endsubhead
Let $\ce_m$ be the set of all $e_*=(e_0,e_1,\do,e_m)\in\NN^{m+1}$ such that $e_0\le e_1\le\do\le e_m$. For any 
$n\in\NN$ let $\ce_m^n=\{e_*\in\ce_m;\sum_ie_i=n\}$.

Let $x_*\in X_m$. We associate to $x_*$ an element $\hx_*\in X_m$ as follows. Let $i_0<i_1<\do<i_s$ be the 
elements of $\fS(x_*)$ in increasing order. Clearly, each of the sets 
$[0,i_0-1],[i_0+1,i_1-1],\do,[i_{s-1}+1,i_s-1],[i_s+1,m]$ has even cardinal, say 
$2t_0,2(t_1-1),\do,2(t_s-1),2t_{s+1}$ (respectively). We define $\hx_*=(\hx_0,\hx_1,\do,\hx_m)\in\NN^{m+1}$ by
$$(\hx_0,\hx_1,\do,\hx_{i_0}-1)=(0,0,1,1,\do,t_0-1,t_0-1), h_{i_0}=t_0,$$
$$\align&(\hx_{i_0+1},\hx_{i_0+2},\do,\hx_{i_1-1})\\&
=(t_0+1,t_0+1,t_0+2,t_0+2,\do,t_0+t_1-1,t_0+t_1-1), h_{i_1}=t_0+t_1,\endalign$$
$$\align&(\hx_{i_1+1},\hx_{i_1+2},\do,\hx_{i_2-1})=\\&
(t_0+t_1+1,t_0+t_1+1,t_0+t_1+2,t_0+t_1+2,\do,t_0+t_1+t_2-1,
t_0+t_1+t_2-1),\\& h_{i_2}=t_0+t_1+t_2,\endalign$$
$$\do$$
$$\align&(\hx_{i_s+1},\hx_{i_s+2},\do,hx_m)\\&
=(t_0+t_1+\do+t_s+1,t_0+t_1+\do+t_s+1,t_0+t_1+\do+t_s+2,\\&
t_0+t_1+\do+t_s+2,\do,t_0+t_1+\do+t_{s+1},t_0+t_1+\do+t_{s+1}).\endalign$$
Note that $\hx_*$ depends only on $\fS(x_*)$, not on $x_*$ itself. We have $\hx_*\in X_m$, $\fS(\hx_*)=\fS(x_*)$.
Let $e_*=x_*-\hx_*$. We have $e_*\in\ce_m$. Moreover for any $i\in[0,m-1]$ such that $\hx_i=\hx_{i+1}$ we have 
$e_i=e_{i+1}$.

\subhead 2.10\endsubhead
Let $x_*\in X_m$, $e_*\in\ce_m$. Then $x_*+e_*\in X_m$ hence $y_*:=x_*+e_*+x_*\in Y_m$. Assume that 
$\fS(e_*+x_*)=\fS(x_*)$ and $(x_*,e_*+x_*)\in S(y_*)$. Then 

$\fS(x_*)\cup\fS(e_*+x_*)=R(y_*)$, $\fS(x_*)\cap\fS(e_*+x_*)=R_0(y_*)$
\nl
hence $\fS(x_*)=R(y_*)=R_0(y_*)$. It follows that for any $\ci\in\fI(y_*)$ we 
have $|\ci|=1$.  

\subhead 2.11\endsubhead
Conversely, let $y_*\in Y_m^n$ be such that for any $\ci\in\fI(y_*)$ we have $|\ci|=1$. By 2.6(a) we can find 
$(x_*,x'_*)\in S(y_*)$. We have $x_*+x'_*=y_*$, $\fS(x_*)\cup\fS(x'_*)=R(y_*)$, $\fS(x_*)\cap\fS(x'_*)=R_0(y_*)$.
From our assumption we have $R_0(y_*)=R(y_*)$. Hence $\fS(x_*)\cup\fS(x'_*)=\fS(x_*)\cap\fS(x'_*)$ so that 
$\fS(x_*)=\fS(x'_*)$. By 2.9 we have $\hx_*=\hx'_*\in X_m$ and $e_*:=x_*-\hx_*\in\ce_m,e'_*:=x'_*-\hx'_*\in\ce_m$.
Moreover, if $i\in[0,m-1]$ and $\hx_i=\hx_{i+1}$ then $e_i=e_{i+1}$ and $e'_i=e'_{i+1}$ hence $\te_i=\te_{i+1}$ 
where $\te_*=e_*+e'_*\in\ce_m$. It follows that $\fS(\te_*+\hx_*)=\fS(\hx_*)=\fS(x_*)=R(y_*)=R_0(y_*)$. Since 
$|\fS(x_*)|+|\fS(x'_*)|=2|\fI(y_*)|$ we have $|\fS(\hx_*)|+|\fS(\te_*+\hx_*)|=2|\fI(y_*)|$. Also, if
$\fI'(y_*)=\em$ then by 2.6(vii) we have $\fS(x'_*)=\em$. Hence $\fS(\hx_*)=\fS(\te_*+\hx_*)=\em$. In any case we
see that $y_*=\hx_*+\te_*+\hx_*$, $(\hx_*,\te_*+\hx_*)\in S(y_*)$,
$\fS(\te_*+\hx_*)=\fS(\hx_*)$.  

\head 3. Type $A_{n-1}$\endhead
\subhead 3.1\endsubhead
For $n\in\NN$ let $S_n$ be the group of all permutations of $\{1,2,\do,n\}$. We have $S_0=S_1=\{1\}$; for $n\ge2$
we regard $S_n$ as a Coxeter group whose generators are the transpositions $(i,i+1)$ for $i\in[1,n-1]$. We have 
$\cs_{S_n}=\Irr(S_n)$. If $k$ is large (relative to $n$) we have a natural bijection $\Irr(S_n)\lra Z_k^n$, 
$[z_*]\lra z_*$, see \cite{\OR, 4.4}. For example, $[(0,1,\do,k-n,k-n+2,\do,k,k+1)]$ is the sign representation
of $S_n$. For any $z_*\in Z_k^n$ we have $\b_0(z_*)=b_{[z_*]}$, see \cite{\OR, (4.4.2)}. 

Assume now that $n=n'+n''$ with $n',n''$ in $\NN$. The set of permutations of $\{1,2,\do,n\}$ which leave stable 
each of the subsets $\{1,2,\do,n'\}$, $\{n'+1,n'+2,\do,n\}$ is a standard parabolic subgroup of $S_n$ which may be
identified with $S_{n'}\T S_{n''}$.

For $z'_*\in Z_k^{n'}$, $z''_*\in Z_k^{n'}$ we have $z'_*+z''_*-z^0_*\in Z_k^n$ and from the definitions we have:

(a) $[[z'_*+z''_*-z^0_*]:\Ind_{S_{n'}\T S_{n''}}^{S_n}([z'_*]\bxt[z''_*])]_{S_n}=1$.
\nl
Note also that $\b_0(z'_*)+\b_0(z''_*)=\b_0(z'_*+z''_*-z^0_*)$ hence
$b_{[z'_*]}+b_{[z''_*]}=b_{[z'_*+z''_*-z^0_*]}$, so that

(b) $[z'_*+z''_*-z^0_*]=j_{S_{n'}\T S_{n''}}^{S_n}([z'_*]\bxt[z''_*])$.

\subhead 3.2\endsubhead
In this subsection we assume that $G$ is of type $A_{n-1}$ ($n\ge2)$. In this case 1.5(a),(b1),(b2) are immediate.
We prove 1.5(b3).

For $C\in\cx$ let $E=\r_C$. We have $E=[z_*]$ for a unique $z_*\in Z_k^n$. We have $\zz_C=1$ and
$\tzz_C=\text{g.c.d.}\{n,z_j-z_j^0(j\in[0,k])\}$ where $\text{g.c.d.}$ denotes the greatest common divisor. We 
identify $\{1,2,\do,n\}=\ZZ/n$ in the obvious way. We also identify $\WW=S_n$ as Coxeter groups so that the
reflections $s_i(i\in\tI)$ are the transpositions $(i,i+1)$ with $i\in\ZZ/n$ (with $i+1$ computed in $\ZZ/n$.) Now
$\Om$ is a cyclic group of order $n$ with generator $\o:i\m i+1$ for all $i\in\ZZ/n$. For any $d|n$ (divisor 
$d\ge1$ of $n$) let $\Om_d$ be the subgroup of $\Om$ generated by $\o^{n/d}$. For any coset $P$ of $\Om_d$ in 
$\Om$ 
let $S_n^P$ be the set of all permutations $w$ of $\ZZ/n$ such that for any $r\in P$ the subset 
$\{r+1,r+2,\do,r+(n/d)\}$ is $w$-stable. We may identify $S_n^P$ with a product of $d$ copies of $S_{n/d}$. Note 
that $\cp^{\Om_d}$ (see 1.11) 
consists of the subgroups $S_n^P$ as above; each of these subgroups is stable under the
conjugation action of $\Om_d$ on $\WW$. An irreducible representation $\bxt_{h=1}^d[\tz_*^{(h)}]$ (with 
$\tz_*^{(h)}\in Z_k^{n/d}$) of $S_n^P$ (identified with $S_{n/d}^d$) is $\Om_d$-stable if and only if 
$\tz_*^{(h)}=\tz_*$ is independent of $h$; in this case we have
$$j_{S_n^P}^{S_n}(\bxt_{h=1}^d[\tz_*^{(h)}])=[\sum_{h=1}^d\tz^{(h)}_*-(d-1)z^0_*]=[d\tz_*-(d-1)z^0_*]$$
as we see by applying $(d-1)$ times 3.1(b). Using this and 1.11 we see that 
$$\fc_{[z_*]}=\max d$$
where $\max$ is taken over all divisors $d\ge1$ of $n$ such that $z_*-z^0_*=d(\tz_*-z^0_*)$ for some 
$\tz_*\in Z_k^{n/d}$. Equivalently, we have
$$\fc_{[z_*]}=\text{g.c.d.}\{n,z_j-z^0_j (j\in[0,k])\}.$$
Since this is equal to $\tzz_C/\zz_C$ we see that 1.5(b3) is proved in our case.

\head 4. Type $B_n$\endhead
\subhead 4.1\endsubhead
For $n\in\NN$ let $W_n$ be the group of permutations of the set \lb
$\{1,2,\do,n,n',\do,2',1'\}$ which commute with 
the involution $i\m i',i'\m i(i\in[1,n])$. We have $W_0=\{1\}$; for $n\ge1$ we regard $W_n$ as a Coxeter group of
type $B_n=C_n$ whose generators are the transposition $(n,n')$ and the products of two transpositions 
$(i,i+1)((i+1)',i')$ for $i\in[1,n-1]$. By \cite{\CL, \S2} we have $\Irr(W_n)=\Irr(W_n)^\da$.

\subhead 4.2\endsubhead
In the remainder of this section we fix an even integer $m=2k$ which is large relative to $n$.

Let $U_k^n=\{(z_*;z'_*)\in Z_k\T Z_{k-1};\r_0(z_*)+\r_0(z'_*)=n\}$. As in \cite{\OR, 4.5} we have a bijection 

(a) $\Irr(W_n)\lra U_k^n$, $[z_*;z'_*]\lra(z_*;z'_*)$.
\nl
(In {\it loc.cit.} the notation $\left(\sm z_*\\z'_*\esm\right)$ was used instead of $(z_*;z'_*)$.) By 
\cite{\CL, \S2} we have 

(b) $b_{[z_*;z'_*]}=2\b_0(z_*)+2\b_0(z'_*)+\r_0(z'_*)$.
\nl
There is a unique bijection $\z_n:\cs_{W_n}@>\si>>X_m^n$ under which $x_*\in X_m^n$ corresponds to$\{[z_*,z'_*]\}$
where $z_*=(x_0,x_2,x_4,\do,x_m)$, $z'_*=(x_1,x_3,x_5,\do,x_{m-1})$. This bijection has the following property: if
$E\in\cs_{W_n},x_*=\z_n(E)$ then $b_E=\b(x_*)$, $f_E=2^{(|\fS(x_*)|-1)/2}$.

\subhead 4.3\endsubhead
Let $u_*\in Z_m$. Define $\ddu_*\in Z_k$, $\du_*\in Z_{k-1}$ by $\ddu_i=u_{2i}-i$ for $i\in[0,k]$, 
$\du_i=u_{2i+1}-i-1$ for $i\in[0,k-1]$. 

\subhead 4.4\endsubhead
Let $(p,q)\in\NN^2$ be such that $p+q=n$. The group of all permutations of $\{1,2,\do,n,n',\do,2',1'\}$ in $W_n$ 
that leave stable each of the subsets
$$\{1,2,\do,p\}, \{p',\do,2',1'\}, \{p+1,\do,n-1,n\}\cup\{n',(n-1)',\do,(p+1)'\}$$
is a standard parabolic subgroup of $W_n$ which may be identified with $S_p\T W_q$ in an obvious way.

Let $(\tz_*;\tz'_*)\in U_k^q$, $u_*\in Z_m^p$. Let 
$v_*=\tz_*+\ddu_*-z^{0,k}_*$, $v'_*=\tz'_*+\du_*-z^{0,k-1}_*$. Then $(v_*;v'_*)\in U_k^n$, $[u_*]\in\Irr(S_p)$, 
$[\tz_*;\tz'_*]\in\Irr(W_q)$, $[v_*;v'_*]\in\Irr(W_n)$. We show:

(a) $[v_*;v'_*]=j_{S_p\T W_q}^{W_n}([u_*]\bxt[\tz_*;\tz'_*])$. 
\nl
We can assume that $p\ge1$ and that the result holds for $p$ replaced by $\tp<p$. In the case where $[u_*]$ is the
sign representation of $S_p$, (a) can be proved along the lines of 
\cite{\UC, 2.7}. If $[u_*]$ is not the sign 
representation of $S_p$, we can find $p',p''$ in $\NN_{>0}$ such that $p'+p''=p$ and $u'\in Z_{2k}^{p'}$, 
$u'\in Z_{2k}^{p''}$ such that $u_*=u'_*+u''_*-z^{0,m}_*$. By 3.1(b), we have 
$[u_*]=j_{S_{p'}\T S_{p''}}^{S_p}([u'_*]\bxt[u''_*])$. Hence 
$$[u_*]\bxt[\tz_*;\tz'_*]=j_{S_{p'}\T S_{p''}\T W_q}^{S_p\T W_q}([u'_*]\bxt[u''_*]\bxt[\tz_*;\tz'_*])$$
and 
$$\align j_{S_p\T W_q}^{W_n}([u_*]\bxt[\tz_*;\tz'_*])&=
j_{S_p\T W_q}^{W_n}j_{S_{p'}\T S_{p''}\T W_q}^{S_p\T W_q}([u'_*]\bxt[u''_*]\bxt[\tz_*;\tz'_*])\\&=
j_{S_{p'}\T S_{p''}\T W_q}^{W_n}([u'_*]\bxt[u''_*]\bxt[\tz_*;\tz'_*])\\&=
j_{S_{p'}\T W_{p''+q}}^{W_n}j_{S_{p'}\T S_{p''}\T W_q}^{S_{p'}\T W_{p''+q}}
([u'_*]\bxt[u''_*]\bxt[\tz_*;\tz'_*])\\&=
j_{S_{p'}\T W_{p''+q}}^{W_n}([u'_*]\bxt[\tz_*+\ddu''_*-z_*^{0,k};\tz'_*+\du''_*-z_*^{0,k-1}])\\&
=[\tz_*+\ddu''_*+\ddu'_*-2z_*^{0,k};\tz'_*+\du''_*+\du'_*-2z_*^{0,k-1}]\\&
=[\tz_*+\ddu_*-z_*^{0,k};\tz'_*+\du_*-z_*^{0,k-1}].\endalign$$
(We have used the induction hypothesis for $p$ replaced by $p'$ or $p''$.) This proves (a).

\subhead 4.5\endsubhead
In the remainder of this section we assume that $G$ has type $B_n$ ($n\ge2$). We identify $\WW=W_n$ as Coxeter 
groups in the standard way. The reflections $s_j(j\in\tI)$ are the transpositions $(n,n')$, $(1,1')$ and the 
products of two transpositions $(i,i+1)(i',(i+1)')$ for $i\in[1,n-1]$. The group $\Om$ has order $2$ with generator
given by the involution $i\m(n+1-i)',i'\m(n+1-i)$ for $i\in[1,n]$.

Let $(r,p,q)\in\NN^3$ be such that $r+p+q=n$. The group of all permutations of \lb
$\{1,2,\do,n,n',\do,2',1'\}$ in 
$W_n$ that leave stable each of the subsets
$$\align&\{1,2,\do,r\}\cup\{r',\do,2',1'\}, \{r+1,r+2,\do,r+p\}, \\&
\{(r+p)',\do,(r+2)',(r+1)'\}, \\&
\{r+p+1,\do,n-1,n\}\cup\{n',(n-1)',\do,(r+p+1)'\}\endalign$$
is a parahoric subgroup of $\WW$ which may be identified with $W_r\T S_p\T W_q$ in an obvious way.

Let $(z_*;z'_*)\in U_k^r$, $(\tz_*;\tz'_*)\in U_k^q$, $u_*\in Z_{2k}^p$, Define $\ddu_*\in Z_k$,$\du_*\in Z_{k-1}$
as in 4.3. Let $w_*=z_*+\tz_*+\ddu_*-2z^{0,k}_*$, $w'_*=z'_*+\tz'_*+\du_*-2z^{0,k-1}_*$. Then 
$(w_*,w'_*)\in U_k^n$, $[z_*;z'_*]\in\Irr(W_r)$, $[u_*]\in\Irr(S_p)$, $[\tz_*;\tz'_*]\in\Irr(W_q)$, 
$[w_*,w'_*]\in\Irr(W_n)$. We show:

(a) $[w_*;w'_*]=j_{W_r\T S_p\T W_q}^{W_n}([z_*;z'_*]\bxt[u_*]\bxt[\tz_*;\tz'_*])$. In particular,

$[[w_*;w'_*]:\Ind_{W_r\T S_p\T W_q}^{W_n}([z_*;z'_*]\bxt[u_*]\bxt[\tz_*;\tz'_*])]_{W_n}=1$.
\nl
Assume first that $p=0$. We have:
$$[[z_*+\tz_*-z^{0,k}_*;z'_*+\tz'_*-z^{0,k-1}_*]:\Ind_{W_r\T W_q}^{W_n}([z_*;z'_*]\bxt[\tz_*;\tz'_*])]_{W_n}=1.
\tag b$$
Using the definitions this can be deduced from the analogous statement for $S_n$, see 3.1(a). Moreover we have 
$b_{[z_*+\tz_*-z^{0,k}_*;z'_*+\tz'_*-z^{0,k-1}_*]}=b_{[z_*;z'_*]}+b_{[\tz_*;\tz'_*]}$. It follows that
$$[z_*+\tz_*-z^{0,k}_*;z'_*+\tz'_*-z^{0,k-1}_*]=j_{W_r\T W_q}^{W_n}([z_*;z'_*]\bxt[\tz_*;\tz'_*]).\tag c$$
Thus (a) holds in this special case. 

In the general case we use 4.4(a) with $n$ replaced by $p+q$ and (c) applied to $n,r,0,p+q$ instead of $n,r,p,q$.
We obtain
$$\align&j_{W_r\T S_p\T W_q}^{W_n}([z_*;z'_*]\bxt[u_*]\bxt[\tz_*;\tz'_*])\\&=
j_{W_r\T W_{p+q}}^{W_n}(j_{W_r\T S_p\T W_q}^{W_r\T W_{p+q}}([z_*;z'_*]\bxt[u_*]\bxt[\tz_*;\tz'_*]))\\&=
j_{W_r\T W_{p+q}}^{W_n}([z_*;z'_*]\bxt[\tz_*+\ddu_*-z^{0,k}_*;\tz'_*+\du_*-z^{0,k-1}_*])=[w_*;w'_*].\endalign$$
This proves (a).

\subhead 4.6\endsubhead
By \cite{\ICC, \S13}, there is a unique bijection $\t:\tcs_\WW@>>>Y_m^n$ such that for any $y_*\in Y_m^n$, the 
fibre $\t\i(y_*)$ is $[z_*,z'_*]$ where $z_*=(y_0,y_2-1,y_4-2,\do,y_m-m/2)$, 
$z'_*=(y_1,y_3-1,y_5-2,\do,y_{m-1}-(m-2)/2)$. This bijection has the following property: if $C\in\cx$ and 
$y_*=\t(\r_C)$, then $\bb_C=\b'(y_*)$, $\zz_C=2^{|\fI(y_*)|-1}$. From \cite{\ICC, \S14} we see that:

$\tzz_C/\zz_C=2$ if $|\ci|=1$ for any $\ci\in\fI(y_*)$,

$\tzz_C/\zz_C=1$ if $|\ci|>1$ for some $\ci\in\fI(y_*)$.

\subhead 4.7\endsubhead
In the setup of 4.5 we assume that $[z_*;z'_*]\in\cs_{W_r}$, $[\tz_*;\tz'_*]\in\cs_{W_q}$. Define $x_*\in X_m^r$,
$\tx_*\in X_m^q$ by $\z_r([z_*;z'_*])=x_*$, $\z_q([\tz_*;\tz'_*])=\tx_*$. Let $e_*=u_*-z^{0,m}_*\in\ce_m$. We 
show: 

(a) $[w_*,w'_*]\in\tcs_\WW$ and $\t([w_*,w'_*])=x_*+e_*+\tx_*$.
\nl
We have $w_i=x_{2i}+\tx_{2i}+u_{2i}-i-2i$ for $i\in[0,k]$, $w'_i=x_{2i+1}+\tx+_{2i+1}+u_{2i+1}-i-1-2i$ for 
$i\in[0,k-1]$. Define $y_*\in\NN^{m+1}$ by $w_i=y_{2i}-i$ for $i\in[0,k]$, $w'_i=y_{2i+1}-i$ for $i\in[0,k-1]$.
Then $y_*=x_*+\tx_*+e_*$. Since $x_*\in X_m,\tx_*\in X_m,e_*\in\ce_m$ we have $y_*\in Y_m$. More precisely,
$y_*\in Y_m^n$. Using 4.6 we deduce that $[w_*,w'_*]\in\tcs_\WW$ and (a) follows.

From (a) and 4.5(a) we see that for $(r,p,q)$ as in 4.5, the assignment 

$(E_1,E_2,\tE_1)\m j_{W_r\T S_p\T W_q}^{W_n}(E_1\bxt E_2\bxt\tE_1)$
\nl
is a map $j:\cs_{W_r}\T\cs_{S_p}\T\cs_{W_q}@>>>\tcs_\WW$ and we have a commutative diagram
$$\CD
\cs_{W_r}\T\cs_{S_p}\T\cs_{W_q}@>j>>\tcs_\WW\\
@V\z_r\T\x_p\T\z_qVV @V\t VV\\
X_m^r\T\ce_m^p\T X_m^q@>h>>Y_m^n \endCD$$                          
where $h$ is given by $(x_*,e_*,\tx_*)\m x_*+e_*+\tx_*$ and $\x_p:\cs_{S_p}@>>>\ce_m^p$ is the bijection 
$[e_*+z_*^{0,m}]\lra e_*$.

\subhead 4.8\endsubhead
Note that $\cp'$ (see 1.9) is exactly the collection of parahoric subgroups $W_r\T S_0\T W_q$ of $W_n$ with 
$(r,p,q)$ as in 4.5 and $p=0$. By 4.7, $j_{W_r\T S_0\T W_q}^{W_n}$ carries $\cs_{W_r}\T\cs_{S_0}\T\cs_{W_q}$ into
$\tcs_\WW$. Hence $\bcs_\WW\sub\tcs_\WW$.

Conversely, let $E\in\tcs_\WW$. With $\t$ as in 4.6, let $y_*=\t(E)\in Y_m^n$. By 2.6(a) we can find 
$(x_*,\tx_*)\in S(y_*)$. Define $r,q$ in $\NN$ by $x_*\in X_m^r,\tx_*\in X_m^q$. We must have $r+q=n$. Let 
$e_*=(0,0,\do,0)\in\ce_m^0$. In the commutative diagram in 4.7 (with $p=0$) we have $h(x_*,e_*,\tx_*)=y_*$, 
$(x_*,e_*,\tx_*)=(\z_r(E_1),\x_p(\QQ),\z_q(\tE_1))$ where $E_1\in\cs_{W_r}$, $\tE_1\in\cs_{W_q}$ (recall that 
$\z_r,\z_q$ are bijections) and $\t(j(E_1,\QQ,\tE_1))=\t(E)$. Since $\t$ is bijective we deduce that 
$E=j(E_1,\QQ,\tE_1)$. Thus, $E\in\bcs_\WW$. Thus, $\tcs_\WW\sub\bcs_\WW$. We see that $\tcs_\WW=\bcs_\WW$. This 
proves 1.5(a) in our case.

\subhead 4.9\endsubhead
In the remainder of this section we fix $C\in\cx$ and we set $E=\r_C\in\tcs_\WW$, $y_*=\t(E)\in Y_m^n$ ($\t$ as in
4.6).

Let $(r,q)\in\NN^2$, $E_1\in\cs_{W_r}$, $\tE_1\in\cs_{W_q}$ be such that $r+q=n$, 

$E=j_{W_r\T S_0\T W_q}^{W_n}(E_1\bxt\QQ\bxt\tE_1)$. 
\nl
(These exist since $E\in\bcs_\WW$.) We set
$x_*=\z_r(E_1)\in X_m^r$, $\tx_*=\z_q(\tE_1)\in X_m^q$. From the commutative diagram in 4.7 we see that 
$x_*+\tx_*=y_*$. By 4.6 we have $\bb_C=\b'(y_*)$. Since $\b'(x_*+\tx_*)=\b(x_*)+\b(\tx_*)$, we have 
$\bb_C=\b(x_*)+\b(\tx_*)$. Since $\b(x_*)=b_{E_1}$, $\b(\tx_*)=b_{\tE_1}$, we have $\bb_C=b_{E_1}+b_{\tE_1}$ hence
$\bb_C=b_{E_1\bxt\QQ\bxt\tE_1}$. Since $E=j_{W_r\T S_0\T W_q}^{W_n}(E_1\bxt\QQ\bxt\tE_1)$ we have 
$b_{E_1\bxt\QQ\bxt\tE_1}=b_E$ hence $\bb_C=b_E$, proving 1.5(b1) in our case. 

Next we note that $f_{E_1}=2^{(|\fS(x_*)|-1)/2}$, $f_{\tE_1}=2^{(|\fS(\tx_*)|-1)/2}$, $\zz_C=2^{|\fI(y_*)|-1}$, 
$|\fS(x_*)|+|\fS(\tx_*)|\le2|\fI(y_*)|$. Hence
$$f_{E_1\bxt\QQ\bxt\tE_1}=2^{(|\fS(x_*)|+|\fS(\tx_*)|-2)/2}\le 2^{|\fI(y_*)|-1}=\zz_C.$$
Taking maximum over all $r,q,E_1,\tE_1$  as above we obtain $\fa_E\le\zz_C$.

Using again 2.6(a) we can find $(x_*,\tx_*)\in S(y_*)$. Define $r,q$ in $\NN$ by $x_*\in X_m^r$, $\tx_*\in X_m^q$.
We must have $r+q=n$. Define $E_1\in\cs_{W_r}$, $\tE_1\in\cs_{W_q}$ by $x_*=\z_r(E_1),\tx_*=\z_q(\tE_1)$. As 
earlier in the proof we have $E=j_{W_r\T\QQ\T W_q}^{W_n}(E_1\bxt\QQ\bxt\tE_1)$. We have 
$$f_{E_1\bxt\QQ\bxt\tE_1}=2^{(|\fS(x_*)|+|\fS(\tx_*)|-2)/2}=2^{|\fI(y_*)|-1}=\zz_C.$$ 
It follows that $\fa_E=\zz_C$, proving 1.5(b2) in our case.

\subhead 4.10\endsubhead
Assume now that $\tzz_C/\zz_C=2$. By 4.6, for any $\ci\in\fI(y_*)$ we have $|\ci|=1$. By 2.11 we can find 
$(r,p,q)$ as in 4.5 with $q=r$ and $x_*\in X_m^r$, $e_*\in\ce_m^p$ such that $y_*=x_*+e_*+x_*$, 
$(x_*,e_*+x_*)\in S(y_*)$, $\fS(e_*+x_*)=\fS(x_*)$. Define $E_1\in\cs_{W_r}$, $E_2\in\cs_{S_p}$ by
$x_*=\z_r(E_1)$, $e_*=\x_p(E_2)$. Using the commutative diagram in 4.7 we see that 
$E=j_{W_r\T S_p\T W_r}^{W_n}(E_1\bxt E_2\bxt E_1)$. Moreover, 
$$f_{E_1\bxt E_2\bxt E_1}=2^{(|\fS(x_*)|+|\fS(x_*)|-2)/2}=2^{(|\fS(x_*)|+|\fS(e_*+x_*)|-2)/2}=2^{|\fI(y_*)|-1}
=\zz_C.$$ 
We have $W_r\T S_p\T W_r=\WW_J$ for a unique $J$ which is $\Om$-stable. Moreover, $E_1\bxt E_2\bxt E_1$ is 
$\Om$-stable. We see that $\fc_E=2$.

\subhead 4.11\endsubhead
Conversely, assume that $\fc_E=2$. Using 1.11 we see that there exist $(r,p,q)$ as in 4.5 with $q=r$ and 
$E_1\in\cs_{W_r}$, $E_2\in\cs_{S_p}$ such that $E=j_{W_r\T S_p\T W_r}^{W_n}(E_1\bxt E_2\bxt E_1)$, 
$f_{E_1\bxt E_2\bxt E_1}=\zz_C$. We set $x_*=\z_r(E_1)\in X_m^r$, $e_*=\x_p(E_2)$. We have $y_*=x_*+e_*+x_*$ and
$$2^{(|\fS(x_*)|+|\fS(x_*)|-2)/2}=2^{|\fI(y_*)|-1};$$
hence $|\fS(x_*)|+|\fS(x_*)|=|\fI(y_*)|$. Let $E'_1=j_{S_p\T W_r}^{W_{p+r}}(E_2\bxt E_1)\in\cs_{W_{p+r}}$. Then 
$E=j_{W_r\T W_{p+r}}^{W_n}(E_1\bxt E'_1)$. Using 1.5(b2) and the definition we have 
$f_{E_1\bxt E'_1}\le\fa_E=\zz_C$. By 1.9(b) we have $f_{E_2\bxt E_1}\le f_{E'_1}$. Hence 
$\zz_C=f_{E_1\bxt E_2\bxt E_1}\le f_{E_1\bxt E'_1}\le\zz_C$; this forces $f_{E_2\bxt E_1}=f_{E'_1}$. The last 
equality can be rewritten as
$$2^{(|\fS(x_*)|-1)/2}=2^{(|\fS(e_*+x_*)|-1)/2}$$
since $e_*+x_*=\z_{p+r}(E'_1)$ (a consequence of 4.4(a)). Hence $|\fS(e_*+x_*)|=|\fS(x_*)|$ and 
$|\fS(x_*)|+|\fS(e_*+x_*)|=2|\fI(y_*)|$. Thus, $(x_*,e_*+x_*)\in S(y_*)$. Using 2.10 we see that for any 
$\ci\in\fI(y_*)$ we have $|\ci|=1$. By 4.6 we have $\tzz_C/\zz_C=2$. 

\subhead 4.12\endsubhead
From 4.10, 4.11, we see that $\tzz_C/\zz_C=2$ if and only if $\fc_E=2$. Since $\fc_E\in[1,2]$ and 
$\tzz_C/\zz_C\in[1,2]$ we see that $\tzz_C/\zz_C=\fc_E$; this proves 1.5(b3) in our case.

\head 5. Type $C_n$\endhead
\subhead 5.1\endsubhead
For $n\in\NN$ let $W'_n$ be the set of all elements in $W_n$ which are even permutations of
$\{1,2,\do,n,n',\do,2',1'\}$. We have $W'_0=W'_1=\{1\}$. For $n\ge2$ we regard $W'_n$ as a Coxeter group of type 
$D_n$ whose generators are the products of two transpositions $(i,i+1)((i+1)',i')$ for $i\in[1,n-1]$ and
$(n-1,n')(n,(n-1)')$.

\subhead 5.2\endsubhead
In this subsection we fix an integer $k$ which is large relative to $n$.

Let $V_k^n$ be the set of unordered pairs $(z_*,z'_*)$ in $Z_{k-1}\T Z_{k-1}$ such that $\r_0(z^*)+\r_0(z'_*)=n$.
If $n\ge2$ we have as in \cite{\OR, 4.5} a map $\io:\Irr(W'_n)@>>>V_k^n$. (In {\it loc.cit.} the notation 
$\left(\sm z_*\\z'_*\esm\right)$ was used instead of $(z_*,z'_*)$.) Now $\io$ is also defined when 
$n\in\{0,1\}$; it is the unique map between two sets of cardinal $1$. 

Let ${}^\da V_k^n$ be the set of {\it ordered} pairs $(z_*;z'_*)$ in $Z_{k-1}\T Z_{k-1}$ such that 
$\r_0(z^*)+\r_0(z'_*)=n$ and either $\r_0(z_*)>\r_0(z'_*)$ or $z_*=z'_*$. We regard ${}^\da V_k^n$ as a subset of
$V_k^n$ by forgetting the order of a pair. We define a partition ${}^\da V_k^n={}'V_k^n\sqc{}''V_k^n$ by

${}''V_k^n=\{(z_*;z'_*)\in{}^\da V_k^n;z_*=z'_*\}$ if $n\ge2$, ${}''V_k^n=\em$ if $n\le1$,

${}'V_k^n=\{(z_*;z'_*)\in{}^\da V_k^n;z_*\ne z'_*\}$ if $n\ge1$, ${}'V_k^n={}^\da V_k^n$ if 
$n=0$.
\nl
By \cite{\CL, \S2} we have $\Irr(W'_n)^\da=\io\i({}^\da V_k^n)$. For $(z_*;z'_*)\in^\da V_k^n$ and $\k\in\{0,1\}$ 
we define $[z_*,z'_*]^\k\in\Irr(W'_n)^\da$ by the following requirements: if $(z_*;z'_*)\in{}'V_k^n$, then 
$\io\i(z_*;z'_*)$ has a single element $[z_*;z'_*]^0=[z_*;z'_*]^1$; if $(z_*;z'_*)\in{}''V_k^n$, then 
$\io\i(z_*;z'_*)$ consists of two elements $[z_*;z'_*]^0,[z_*;z'_*]^1$.

By \cite{\CL, \S2}, if $(z_*;z'_*)\in{}^\da V_k^n$ then $b_{[z_*;z'_*]^\k}=2\b_0(z_*)+2\b_0(z'_*)+\r_0(z'_*)$. 

There is a unique map $\z'_n:\cs_{W'_n}@>>>X_{2k-1}^n$ such that for any $x_*\in X_{2k-1}^n$, $\z'_n{}\i(x_*)$ is
$\{[z_*;z'_*]^0=[z_*;z'_*]^1\}$ (if $\fS(x_*)\ne\em$ or if $n=0$) and is $\{[z_*;z'_*]^0,[z_*;z'_*]^1\}$ (if
$\fS(x_*)=\em$ and $n\ge2$) where 

$z_*=(x_1,x_3,x_5,\do,x_{2k-1})$, $z'_*=(x_0,x_2,x_4,\do,x_{2k-2})$. 
\nl
This map has the following property: if $E\in\cs_{W'_n},x_*=\z'_n(E)$ then $b_E=\b(x_*)$, 
$f_E=2^{\max((|\fS(x_*)|-2)/2,0)}$.

There is a unique map $\tiz_n:\cs_{W'_n}@>>>\tX_{2k}^n$ such that for any $x_*\in\tX_{2k}^n$, $\tiz_n\i(x_*)$ is
$\{[z_*;z'_*]^0=[z_*,z'_*]^1\}$ (if $\fS(x_*)\ne\{0\}$ or if $n=0$) and is $\{[z_*;z'_*]^0,[z_*,z'_*]^1\}$ (if 
$\fS(x_*)=\{0\}$ and $n\ge2$) where $z_*=(x_2-1,x_4-1,\do,x_{2k}-1)$, $z'_*=(x_1-1,x_3-1,x_5-1,\do,x_{2k-1}-1)$.

This map has the following property: if $E\in\cs_{W'_n},x_*=\tiz_n(E)$, then $b_E=\b(x_*)$, 
$f_E=2^{\max((|\fS(x_*)|-3)/2,0)}$.

\subhead 5.3\endsubhead
In the remainder of this section we assume that $G$ is of type $C_n$ ($n\ge3$) and we identify $\WW=W_n$ as 
Coxeter groups in the standard way; we also fix an even integer $m=2k$ which is large relative to $n$. The 
reflections $s_j(j\in\tI)$ are the transposition $(1,1')$ and the products of two transpositions 
$(i,i+1)(i',(i+1)')$ for $i\in[1,n-1]$ and $(n-1,n')(n,(n-1)')$. The group $\Om$ has order $2$ with generator given
by the transposition $(n,n')$.

Let $(r,q)\in\NN^2$ be such that $r+q=n$. The group of all permutations of $\{1,2,\do,n,n',\do,2',1'\}$ in $W_n$ 
that leave stable the subset $\{1,2,\do,r\}\cup\{r',\do,2',1'\}$ and which restrict to an even permutation of 
$\{r+1,\do,n-1,n\}\cup\{n',(n-1)',\do,(r+1)'\}$, is a parahoric subgroup of $\WW$ which may be identified with 
$W_r\T W'_q$ in an obvious way. Let $(z_*;z'_*)\in U_k^r$, 
$(\tz_*;\tz'_*)\in{}^\da V_k^q$. Let 

$\tz_*^!=(0,\tz_0+1,\tz_1+1,\do,\tz_{k-1}+1)\in Z_k$. 
\nl
Let $w_*=z_*+\tz_*^!-z^{0,k}_*$, $w'_*=z'_*+\tz'_*-z^{0,k-1}_*$. Then $[z_*;z'_*]\in\Irr(W_r)$, 
$[\tz_*;\tz'_*]^\k\in\Irr(W'_q)^\da$ $(k=0,1$), $[w_*;w'_*]\in\Irr(W_n)$ are well defined and we have

(a) $[[w_*;w'_*]:\Ind_{W_r\T W'_q}^{W_n}([z_*;z'_*]\bxt[\tz_*;\tz'_*]^\k)]_{W_n}=1$.
\nl
(This can be deduced from the second sentence in 4.5(a) with $p=0$.) Moreover, we have 
$b_{[w_*;w'_*]}=b_{[z_*;z'_*]}+b_{[\tz_*;\tz'_*]^\k}$. It follows that

(b) $[w_*;w'_*]=j_{W_r\T W'_q}^{W_n}([z_*;z'_*]\bxt[\tz_*;\tz'_*]^\k)$.

\subhead 5.4\endsubhead
By \cite{\ICC, \S12}, there is a unique bijection $\ti\t:\tcs_\WW@>>>\tY_m^n$ such that for any $y_*\in\tY_m^n$, 
the fibre $\ti\t\i(y_*)$ is $\{[z_*,z'_*]\}$ where $z_*=(y_0,y_2-1,y_4-2,\do,y_m-m/2)$, 
$z'_*=(y_1-1,y_3-2,y_5-3,\do,y_{m-1}-m/2)$. This bijection has the following property: if $C\in\cx$ and 
$y_*=\ti\t(\r_C)$ then $\bb_C=\tib'(y_*)$, $\zz_C=2^{|\fI(y_*)|-1-\ti\d_{y_*})}$ where $\ti\d_{y_*}=1$ if there 
exists $\ci\in\fI'(y_*)$ such that $0\n\ci$ and $\ti\d_{y_*}=0$ if there is no $\ci\in\fI'(y_*)$ such that 
$0\n\ci$. Moreover, $\tzz_C=2^{|\fI(y_*)|-1)}$. Hence $\tzz_C/\zz_C=2^{\ti\d_{y_*}}$.

\subhead 5.5\endsubhead
In the setup of 5.3 we assume that $[z_*;z'_*]\in\cs_{W_r}$, $[\tz_*;\tz'_*]^k\in\cs_{W'_q}$. We set 
$x_*=\z_r([z_*;z'_*])\in X_m^r$, $\tx_*=\tiz_q([\tz_*;\tz'_*]^\k)\in\tX_m^q$. We show:

(a) $[w_*,w'_*]\in\tcs_\WW$ and $\ti\t([w_*,w'_*])=x_*+\tx_*$.
\nl
We have $w_i=x_{2i}+\tx_{2i}-i$, $w'_i=x_{2i+1}+\tx_{2i+1}-1-i$. Define $y_*\in\NN^{m+1}$ by $y_{2i}=w_i+i$ for 
$i\in[0,k]$, $y_{2i+1}=w'_i+i+1$ for $i\in[0,k-1]$. We have $y_*=x_*+\tx_*$. Since $x_*\in X_m,\tx_*\in X_m,$ we 
have $y_*\in\tY_m$. More precisely, $y_*\in\tY_m^n$. Using 5.4 we deduce that $[w_*,w'_*]\in\tcs_\WW$ and (a) 
follows.

From (a) and 5.3(b) we see that for $(r,q)$ as in 5.3, the assignment 
$(E_1,\tE_1)\m j_{W_r\T W'_q}^{W_n}(E_1\bxt\tE_1)$ is a map $j:\cs_{W_r}\T\cs_{W'_q}@>>>\tcs_\WW$ and we have a 
commutative diagram
$$\CD
\cs_{W_r}\T\cs_{W'_q}@>j>>\tcs_\WW\\
@V\z_r\T\tiz_qVV @V\ti\t VV\\
X_m^r\T\tX_m^q@>h>>\tY_m^n\endCD$$                          
where $h$ is given by $(x_*,\tx_*)\m x_*+\tx_*$.

\subhead 5.6\endsubhead
Note that $\cp'$ is exactly the collection of subgroups $W_r\T W'_q$ of $W_n$ with $(r,q)$ as in 5.3 and $q\ne1$.
(On the other hand $W_{n-1}\T W'_1$ is a maximal parabolic subgroup of the Coxeter group $W_n$.) By 5.5, 
$j_{W_r\T W'_q}^{W_n}$ carries $\cs_{W_r}\T\cs_{W'_q}$ into $\tcs_\WW$. Hence $\bcs_\WW\sub\tcs_\WW$.

Conversely, let $E\in\tcs_\WW$. With $\ti\t$ as in 5.4, let $y_*=\ti\t(E)\in\tY_m^n$. By 2.7(b) we can find 
$(x_*,\tx_*)\in\tS(y_*)$. (The assumption 2.7(a) is automatically satisfied since $m$ is large relative to $n$.) 
Define $r,q$ in $\NN$ by $x_*\in X_m^r,\tx_*\in\tX_m^q$. We must have $r+q=n$. In the commutative diagram in 5.5 
we have $h(x_*,\tx_*)=y_*$, $(x_*,\tx_*)=(\z_r(E_1),\tiz_q(\tE_1))$ where $E_1\in\cs_{W_r}$, $\tE_1\in\cs_{W'_q}$
(recall that $\z_r$, $\tiz_q$ are surjective) and $\ti\t(j(E_1,\tE_1))=\ti\t(E)$. Since $\ti\t$ is bijective we 
deduce that $E=j(E_1,\tE_1)$. Thus, $E\in\bcs_\WW$ and $\tcs_\WW\sub\bcs_\WW$. We see that $\tcs_\WW=\bcs_\WW$. 
This proves 1.5(a) in our case.

\subhead 5.7\endsubhead
In the remainder of this section we fix $C\in\cx$ and we set $E=\r_C\in\tcs_\WW$, $y_*=\ti\t(E)\in\tY_m^n$ (with 
$\ti\t$ as in 5.4). 

Let $(r,q)\in\NN^2$, $E_1\in\cs_{W_r}$, $\tE_1\in\cs_{W'_q}$ be such that $r+q=n$, 
$E=j_{W_r\T W'_q}^{W_n}(E_1\bxt\tE_1)$. (These exist since $E\in\bcs_\WW$.) We set $x_*=\z_r(E_1)\in X_m^r$,
$\tx_*=\tiz_q(\tE_1)\in\tX_m^q$. From the commutative diagram in 5.5 we see that $x_*+\tx_*=y_*$. By 5.4 we have 
$\bb_C=\tib'(y_*)$. Since $\tib'(x_*+\tx_*)=\b(x_*)+\tib(\tx_*)$ we have $\bb_C=\b(x_*)+\tib(\tx_*)$. Since 
$\b(x_*)=b_{E_1}$, $\tib(\tx_*)=b_{\tE_1}$, we have $\bb_C=b_{E_1}+b_{\tE_1}$ hence $\bb_C=b_{E_1\bxt\tE_1}$. 
Since $E=j_{W_r\T W'_q}^{W_n}(E_1\bxt\tE_1)$ we have $b_{E_1\bxt\tE_1}=b_E$ hence $\bb_C=b_E$, proving 1.5(b1) in
our case. 

If $|\fS(\tx_*)|\ge3$ then 
$$\align&f_{E_1\bxt\tE_1}=f_{E_1}f_{\tE_1}=2^{(|\fS(x_*)|-1)/2}2^{(|\fS(\tx_*)|-3)/2}\\&=
2^{(|\fS(x_*)|+|\fS(\tx_*)|-4)/2}\le 2^{|\fI(y_*)|-2}\le\zz_C.\endalign$$
If $|\fS(\tx_*)|=1$ and $|\fS(x_*)|\le2|\fI(y_*)|-3$ then 
$$f_{E_1\bxt\tE_1}=f_{E_1}f_{\tE_1}=2^{(|\fS(x_*)|-1)/2}\le 2^{|\fI(y_*)|-2}\le\zz_C.$$
If $|\fS(\tx_*)|=1$ (hence $\fS(\tx_*)=\{0\}$) and $|\fS(x_*)|=2|\fI(y_*)|-1$ then $\ti\d_{y_*}=0$ so that 
$\zz_C=2^{|\fI(y_*)|-1}$ and 
$$f_{E_1\bxt\tE_1}=f_{E_1}f_{\tE_1}=2^{(|\fS(x_*)|-1)/2}=2^{|\fI(y_*)|-1}=\zz_C.$$
Thus in any case we have $f_{E_1\bxt\tE_1}\le\zz_C$. Taking maximum over all $r,q,E_1,\tE_1$ as above we obtain 
$\fa_E\le\zz_C$.

\subhead 5.8\endsubhead
Assume now that $\ti\d_{y_*}=1$. Then $|\fI(y_*)|\ge2$. By 2.7(b) we can find $(x_*,\tx_*)\in\tS(y_*)$. By 2.7(d)
we have $|\fS(\tx_*)|\ge3$. Define $(r,q)\in\NN^2$ by $x_*\in X_m^r,\tx_*\in\tX_m^q$. We must have $r+q=n$. We can
find $E_1\in\cs_{W_r}$, $\tE_1\in\cs_{W'_q}$ such that $x_*=\z_r(E_1),\tx_*=\tiz_q(\tE_1)$. As earlier in the 
proof, we have $E=j_{W_r\T W'_q}^{W_n}(E_1\bxt\tE_1)$ and 
$$f_{E_1\bxt\tE_1}=f_{E_1}f_{\tE_1}=2^{(|\fS(x_*)|-1)/2+(|\fS(\tx_*)|-3)/2}=2^{|\fI(y_*)|-2}=\zz_C.$$

\subhead 5.9\endsubhead
Next we assume that $\ti\d_{y_*}=0$. By 2.7(b) we can find $(x_*,\tx_*)\in\tS(y_*)$. By 2.7(v) we have
$\fS(\tx_*)=\{0\}$. Then $|\fS(x_*)|=2|\fI(y_*)|-1$. Define $(r,q)\in\NN^2$ by $x_*\in X_m^r,\tx_*\in\tX_m^q$. We
must have $r+q=n$. We can find $E_1\in\cs_{W_r}$, $\tE_1\in\cs_{W'_q}$ such that $x_*=\z_r(E_1)$,
$\tx_*=\tiz_q(\tE_1)$. We have $E=j_{W_r\T W'_q}^{W_n}(E_1\bxt\tE_1)$ and 
$$f_{E_1\bxt\tE_1}=f_{E_1}f_{\tE_1}=2^{(|\fS(x_*)|-1)/2}=2^{|\fI(y_*)|-1}=\zz_C.$$
Using this and 5.8 we see that in any case, $\fa_E=\zz_C$, proving 1.5(b2) in our case.

\subhead 5.10\endsubhead
Assume first that $\d_{y_*}=1$. Let $r,q,x_*,\tx_*,E_1,\tE_1$ be as in 5.8. Then $|\fS(\tx_*)|\ge3$ hence $q\ge1$
(so that the unique $J$ such that $\WW_J=W_r\T W'_q$ is $\Om$-stable) and $\tE_1$ is $\Om$-stable. It follows that 
$\fc_E=2$.

Conversely, assume that $\fc_E=2$. Using 1.11 we see that there exist $(r,q)\in\NN^2$ be such that $r+q=n$ with 
$q\ge1$ and $E_1\in\cs_{W_r}$, $\tE_1\in\cs_{W'_q}$ such that $\tE_1$ is $\Om$-stable, 
$E=j_{W_r\T W'_q}^{W_n}(E_1\bxt\tE_1)$, $f_{E_1\bxt\tE_1}=\zz_C$. We set $x_*=\z_r(E_1)\in X_m^r$, 
$\tx_*=\tiz_q(\tE_1)\in\tX_m^q$. We have $y_*=x_*+\tx_*$. Since $\tE_1$ is $\Om$-stable, we have $|\fS(\tx_*)|\ge3$.
Hence 
$$2^{|\fI(y_*)|-2}\le\zz_C=f_{E_1\bxt\tE_1}=2^{(|\fS(x_*)|-1)/2+(|\fS(\tx_*)|-3)/2}\le2^{|\fI(y_*)|-2}.$$
It follows that $2^{|\fI(y_*)|-2}=\zz_C$ so that $\ti\d_{y_*}=1$. 

We see that $\tzz_C/\zz_C=2$ if and only if $\fc_E=2$. Since $\fc_E\in[1,2]$ and $\tzz_C/\zz_C\in[1,2]$ we see 
that $\tzz_C/\zz_C=\fc_E$; this proves 1.5(b3) in our case.
 
\head 6. Type $D_n$\endhead
\subhead 6.1\endsubhead
In this section we assume that $G$ is of type $D_n$ ($n\ge4$). We identify $\WW=W'_n$ as Coxeter groups in the 
usual way. The reflections $s_j(j\in\tI)$ are the products of two transpositions $(i,i+1)(i',(i+1)')$ for 
$i\in[1,n-1]$ and $(n-1,n')(n,(n-1)')$, $(1,2')(2,1')$. Define $\o_1\in W'_n$ by $i\m(n+1-i)',i'\m n+1-i$ for 
$i\in[1,n-1]$, $n\m1$, $n'\m1'$ (if $n$ is even) and by $i\m(n+1-i)',i'\m n+1-i$ for $i\in[1,n]$ (if $n$ is 
even). Define $\o_2\in W'_n$ by $i\m i$ for $i\in[2,n-1]$, $1\m1',1'\m1,n\m n',n'\m n$. We have $\o_1,\o_2\in\Om$.
If $n$ is odd, $\Om$ is cyclic of order $4$ with generator $\o_1$ such that $\o_1^2=\o_2$. If $n$ is even, $\Om$ is
noncyclic of order $4$ with generators $\o_1,\o_2$ of order $2$.

\subhead 6.2\endsubhead
In the remainder of this section we fix an odd integer $m=2k-1$ which is large relative to $n$.

Let $(p,q)\in\NN^2$ be such that $p+q=n$, $q\ge1$. The group of all permutations of $\{1,2,\do,n,n',\do,2',1'\}$ 
in $W_n$ that leave stable each of the subsets $\{1,2,\do,p\}$, $\{p',\do,2',1'\}$ and induce an even permutation
on the subset $\{p+1,\do,n-1,n\}\cup\{n',(n-1)',\do,(p+1)'\}$ is a standard parabolic subgroup of $W'_n$ which may
be identified with $S_p\T W'_q$ in an obvious way.

Let $(\tz_*;\tz'_*)\in{}'V_k^q$, $u_*\in Z_{2k-1}^p$. Define $\ddu_*\in Z_{k-1}$, $\du_*\in Z_{k-1}$ by 
$\ddu_i=u_{2i}-i$, $\du_i=u_{2i+1}-i-1$ for $i\in[0,k-1]$. Let $v_*=\tz_*+\du_*-z^{0,k-1}_*$, 
$v'_*=\tz'_*+\ddu_*-z^{0,k-1}_*$. Then $(v_*;v'_*)\in{}'V_k^n$, $[u_*]\in\Irr(S_p)$, 
$[\tz_*;\tz'_*]\in\Irr(W'_q)$, $[v_*;v'_*]\in\Irr(W'_n)$. We have:

(a) $[v_*;v'_*]^0=j_{S_p\T W'_q}^{W'_n}([u_*]\bxt[\tz_*;\tz'_*]^0)$. 
\nl
The proof is similar to that of 4.4(a).

\subhead 6.3\endsubhead
Let $(r,p,q)\in\NN^3$ be such that $r+p+q=n$. The group of all permutations of $\{1,2,\do,n,n',\do,2',1'\}$ in 
$W'_n$ that leave stable each of the subsets 
$$\{r+1,r+2,\do,r+p\}, \{(r+p)',\do,(r+2)',(r+1)'\}$$
and induce an even permutation on each of the subsets 
$$\{1,2,\do,r\}\cup\{r',\do,2',1'\}, \{r+p+1,\do,n-1,n\}\cup\{n',(n-1)',\do,(r+p+1)'\}$$
is a parahoric subgroup of $\WW$ which may be identified with $W'_r\T S^{(0)}_p\T W'_q$ in an obvious way.
($S^{(0)}_p$ is a copy of $S_p$.)

When $r=0,p\ge2$, the group of all permutations of $\{1,2,\do,n,n',\do,2',1'\}$ in $W'_n$ that leave stable each 
of the subsets 
$$\{1',2,\do,p\},\{p',\do,2',1\},\{p+1,\do,n-1,n\}\cup\{n',(n-1)',\do,(p+1)'\}$$
is a parahoric subgroup of $\WW$ which may be identified with $W'_r\T S_p^{(1)}\T W'_q$. ($S^{(1)}_p$ is a copy of
$S_p$.)

When $p\ge2,q=0$, the group of all permutations of $\{1,2,\do,n,n',\do,2',1'\}$ in $W'_n$ that leave stable each 
of the subsets 
$$\{r+1,r+2,\do,n-1,n'\},\{n,(n-1)'\do,(r+2)',(r+1)'\},\{1,2,\do,r\}\cup\{r',\do,2',1'\}$$ 
is a parahoric subgroup of $\WW$ which may be identified with $W'_r\T S_p^{(2)}\T W'_q$. ($S^{(2)}_p$ is a copy of
$S_p$.)

When $r=q=0$, the group of all permutations of $\{1,2,\do,n,n',\do,2',1'\}$ in $W'_n$ that leave stable each of 
the subsets $\{1',2,3,\do,n-1,n'\}$, $\{n,(n-1)'\do,3',2',1\}$ is a parahoric subgroup of $\WW$ which may be 
identified with $W'_r\T S_p^{(3)}\T W'_q$. ($S^{(3)}_p$ is a copy of $S_p$.)

Thus the parahoric subgroup $W'_r\T S_p^{(\l)}\T W'_q$ is defined in the following cases:

(a) $\l=0$; $p\ge2,r=0,\l=1$; $p\ge2,q=0,\l=2$; $r=q=0,\l=3$.
\nl
When $p=0$ we write also $W'_r\T W'_q$ instead of $W'_r\T S_p^{(0)}\T W'_q$.

Let $(z_*;z'_*)\in{}^\da V_k^r$, $(\tz_*;\tz'_*)\in{}^\da V_k^q$, $u_*\in Z_{2k-1}^p$, Define $\ddu_*\in Z_{k-1}$,
$\du_*\in Z_{k-1}$ by $\ddu_i=u_{2i}-i$, $\du_i=u_{2i+1}-i-1$ for $i\in[0,k-1]$. Let 
$w_*=z_*+\tz_*+\du_*-2z^{0,k-1}_*$, $w'_*=z'_*+\tz'_*+\ddu_*-2z^{0,k-1}_*$. Then $(w_*,w'_*)\in{}^\da V_k^n$.

For $\k,\ti\k,\k'\in\{0,1\}$ we have $[z_*;z'_*]^\k\in\Irr(W'_r)^\da$, $[u_*]\in\Irr(S_p)$, 
$[\tz_*;\tz'_*]^{\ti\k}\in\Irr(W'_q)^\da$, $[w_*,w'_*]^{\k'}\in\Irr(W'_n)^\da$. For $\l$ as in (a) we have:

(b) $[w_*;w'_*]^{\k'}=j_{W'_r\T S_p^{(\l)}\T W'_q}^{W'_n}([z_*;z'_*]^\k\bxt[u_*]\bxt[\tz_*;\tz'_*]^{\ti\k})$
\nl
with the following restriction on $\k'$: if $z_*=z'_*,\tz_*=\tz'_*$, 
$\du_*=\ddu_*$, then $w_*=w'_*$ and $\k'$ in 
(b) is uniquely determined by $\k,\ti\k,\l$; moreover, both $\k'=0$ and $\k'=1$ are obtained from some 
$(\k,\ti\k,\l)$.

Now (b) can be proved in a way similar to 4.5(a); alternatively, from the second statement of 4.5(a) one can 
deduce that
$$[[w_*;w'_*]^{\k'}:
\Ind_{W'_r\T S_p^{(\l)}\T W'_q}^{W'_n}([z_*;z'_*]^\k\bxt[u_*]\bxt[\tz_*;\tz'_*]^{\ti\k})]_{W'_n}\ge1;$$
we can also check directly that $b_{[w_*;w'_*]^{\k'}}=b_{[z_*;z'_*]^\k}+b_{[u_*]}+b_{[\tz_*;\tz'_*]^{\ti\k}}$ 
and (b) follows.

\subhead 6.4\endsubhead
By \cite{\ICC, \S13}, we have $\tcs_\WW\sub\Irr(W'_n)^\da$ and there is a unique map $\t:\tcs_\WW@>>>Y_m^n$ such 
that for $y_*\in Y_m^n$, $\t\i(y_*)$ consists of $[z_*;z'_*]^0=[z_*;z'_*]^1$ (if $\fI(y_*)\ne\em$) and consists of
$[z_*;z'_*]^0,[z_*;z'_*]^1$ (if $\fI(y_*)=\em$) where 

$z_*=(y_1,y_3-1,y_5-2,\do,y_m-(m-1)/2)$, 

$z'_*=(y_0,y_2-1,y_4-2,\do,y_{m-1}-(m-1)/2)$.
\nl
This map has the following property: if $C\in\cx$ and 
$y_*=\t(\r_C)$, then $\bb_C=\b'(y_*)$, $\zz_C=2^{\max(|\fI(y_*)|-1-\d_{y_*},0)}$ where $\d_{y_*}=1$ if
$\fI'(y_*)\ne\em$ and $\d_{y_*}=0$ if $\fI'(y_*)=\em$. Moreover, $\tzz_C/\zz_C$ is:

$4$ if $\d_{y_*}=1$ and $|\ci|=1$ for any $\ci\in\fI(y_*)$,

$2$ if $\d_{y_*}=1$ and $|\ci|>1$ for some $\ci\in\fI(y_*)$,

$2$ if $\fI(y_*)=\em$,

$1$ if $\d_{y_*}=0$ and $|\ci|>1$ for some $\ci\in\fI(y_*)$.
\nl
More precisely, let $\uG@>>>G$ be a double covering which is a special orthogonal group and let $\uzz_C$ be the 
number of connected components of the centralizer in $\uG$ of a unipotent element of $\uG$ which maps to an 
element of $C$. From \cite{\ICC, \S14} we see that:

$\tzz_C/\uzz_C=2$ if $|\ci|=1$ for any $\ci\in\fI(y_*)$,

$\tzz_C/\zz_C=1$ if $|\ci|>1$ for some $\ci\in\fI(y_*)$.
\nl
On the other hand, from $\uzz_C=2^{\max(|\fI(y_*)|-1,0)}$, $\zz_C=2^{\max(|\fI(y_*)|-1-\d_{y_*},0)}$, we see that
$\uzz_C/\zz_C=2^{\d_{y_*}}$.

\subhead 6.5\endsubhead
In the setup of 6.3 we assume that $[z_*;z'_*]^\k\in\cs_{W'_r}$, $[\tz_*;\tz'_*]^{\ti\k}\in\cs_{W'_q}$ and $\k'$ 
is as in 6.3(b). Define $x_*\in X_m^r,\tx_*\in X_m^q$ by $\z'_r([z_*;z'_*]^\k)=x_*$, 
$\z'_q([\tz_*;\tz'_*]^{\ti\k})=\tx_*$. Let $e_*=u_*-z^{0,m}_*\in\ce_m$. We show: 

(a) $[w_*,w'_*]^{\k'}\in\tcs_\WW$ and $\t([w_*,w'_*]^{\k'})=x_*+e_*+\tx_*$.
\nl
We have $w_i=x_{2i+1}+\tx_{2i+1}+u_{2i+1}-i-2i-1$, $w'_i=x_{2i}+\tx_{2i}+u_{2i}-i-2i$ for $i\in[0,k-1]$. Define 
$y_*\in\NN^{m+1}$ by $w_i=y_{2i+1}-i$, $w'_i=y_{2i}-i$ for $i\in[0,k-1]$. Then $y_*=x_*+\tx_*+e_*$. Since 
$x_*\in X_m,\tx_*\in X_m,e_*\in\ce_m$ we have $y_*\in Y_m$. More precisely, $y_*\in Y_m^n$. Using 6.4 we deduce
that $[w_*,w'_*]^{\k'}\in\tcs_\WW$ and (a) follows.

From (a) and 6.3(b) we see that for $\l$ as in 6.3(a), the assignment 
$(E_1,E_2,\tE_1)\m j_{W_r\T S_p^{(\l)}\T W_q}^{W_n}(E_1\bxt E_2\bxt\tE_1)$ is a map
$j:\cs_{W'_r}\T\cs_{S_p}\T\cs_{W'_q}@>>>\tcs_\WW$ and we have a commutative diagram
$$\CD
\cs_{W'_r}\T\cs_{S_p}\T\cs_{W'_q}@>j>>\tcs_\WW\\
@V\z'_r\T\x_p\T\z'_qVV @V\t VV\\
X_m^r\T\ce_m^p\T X_m^q@>h>>Y_m^n \endCD$$                          
where $h$ is given by $(x_*,e_*,\tx_*)\m x_*+e_*+\tx_*$ and $\x_p:\cs_{S_p}@>>>\ce_m^p$ is the bijection 
$[e_*+z_*^{0,m}]\lra e_*$.

\subhead 6.6\endsubhead
Note that $\cp'$ is exactly the collection of parahoric subgroups $W'_r\T W'_q$ of $W'_n$ with 
$(r,q)\in\NN^2$ such that $r+q=n$, $r\ne1$, $q\ne1$. (On the other hand $W'_{n-1}\T W'_1$, $W'_1\T W'_{n-1}$ are maximal 
parabolic subgroup of the Coxeter group $W'_n$.) By 6.5, $j_{W'_r\T W'_q}^{W'_n}$ carries $\cs_{W'_r}\T\cs'_{W_q}$
into $\tcs_\WW$. Hence $\bcs_\WW\sub\tcs_\WW$.

Conversely, let $E\in\tcs_\WW$, $y_*=\t(E)\in Y_m^n$ ($\t$ as in 6.4). By 2.6(a) we can find 
$(x_*,\tx_*)\in S(y_*)$. Define $r,q$ in $\NN$ by $x_*\in X_m^r$, $\tx_*\in X_m^q$. We must have $r+q=n$. In the 
commutative diagram in 6.5 we have $h(x_*,\tx_*)=y_*$, $x_*=\z'_r(E_1),\tx_*=\z'_q(\tE_1)$ where 
$E_1\in\cs_{W'_r}$, $\tE_1\in\cs_{W'_q}$ (recall that $\z'_r$, $\z'_q$ are surjective) and
$\t(j(E_1,\tE_1)=\t(E)$. Thus $j(E_1,\tE_1),E$ are in the same fibre of $\io:{}^\da V_k^n@>>>\Irr(W'_n)^\da$. 
Replacing $E_1$ or $\tE_1$ by an element in the same fibre of $\io:{}^\da V_k^r@>>>\Irr(W'_r)^\da$ or 
$\io:{}^\da V_k^q@>>>\Irr(W'_q)^\da$ we see that we can assume that $j(E_1,\tE_1)=E$. Thus, $E\in\bcs_\WW$. Thus, $\tcs_\WW\sub\bcs_\WW$. We see that $\tcs_\WW=\bcs_\WW$. This proves 1.5(a) in our case.

\subhead 6.7\endsubhead
In the remainder of this section we fix $C\in\cx$ and we set $E=\r_C\in\tcs_\WW$, $y_*=\t(E)\in Y_m^n$ (with
$\t$ as in 6.4).

Let $(r,q)\in\NN^2$, $E_1\in\cs_{W'_r}$, $\tE_1\in\cs_{W'_q}$ be such that $r+q=n$, 
$E=j_{W'_r\T W'_q}^{W'_n}(E_1\bxt\tE_1)$. (These exist since $E\in\bcs_\WW$.) Define $x_*\in X_m^r,\tx_*\in X_m^q$
by $x_*=\z'_r(E_1),\tx_*=\z'_q(\tE_1)$. From the commutative diagram in 6.5 we see that $x_*+\tx_*=y_*$. By 6.4,
we have $\bb_C=\b'(y_*)$. Since $\b'(x_*+\tx_*)=\b(x_*)+\b(\tx_*)$ we have $\bb_C=\b(x_*)+\b(\tx_*)$. Since
$\b(x_*)=b_{E_1}$, $\b(\tx_*)=b_{\tE_1}$, we have $\bb_C=b_{E_1}+b_{\tE_1}$ hence $\bb_C=b_{E_1\bxt\tE_1}$. Since 
$E=j_{W'_r\T W'_q}^{W'_n}(E_1\bxt\tE_1)$ we have $b_{E_1\bxt\tE_1}=b_E$ hence $\bb_C=b_E$, proving 1.5(b1) in our
case.

If $|\fS(x_*)|\ge2$, $|\fS(\tx_*)|\ge2$, then
$$\align&f_{E_1\bxt\tE_1}=f_{E_1}f_{\tE_1}=2^{(|\fS(x_*)|-2)/2}2^{(|\fS(\tx_*)|-2)/2)}\\&=
2^{(|\fS(x_*)|+|\fS(\tx_*)|-4)/2}\le 2^{|\fI(y_*)|-2}\le\zz_C.\endalign$$ 
If $|\fS(x_*)|=0$, $2\le|\fS(\tx_*)|\le2|\fI(y_*)|-2$ then
$$f_{E_1\bxt\tE_1}=f_{E_1}f_{\tE_1}=2^{(|\fS(\tx_*)|-2)/2)}\le 2^{|\fI(y_*)|-2}\le\zz_C.$$ 
Similarly, if $2\le|\fS(x_*)|\le2|\fI(y_*)|-2$, $|\fS(\tx_*)|=0$, then $f_{E_1\bxt\tE_1}\le\zz_C$. If 
$|\fS(x_*)|=0$, $2\le|\fS(\tx_*)|=2|\fI(y_*)|$ then $\fI'(y_*)=\em$ and 
$$f_{E_1\bxt\tE_1}=f_{E_1}f_{\tE_1}=2^{(|\fS(\tx_*)|-2)/2)}=2^{|\fI(y_*)|-1}\le\zz_C.$$ 
Similarly, if $2\le|\fS(x_*)|=2|\fI(y_*)|$, $|\fS(\tx_*)|=0$ then $f_{E_1\bxt\tE_1}\le\zz_C$. If $|\fS(x_*)|=0$, 
$|\fS(\tx_*)|=0$, then 
$$f_{E_1\bxt\tE_1}=f_{E_1}f_{\tE_1}=1=\zz_C.$$ 
Thus in any case we have $f_{E_1\bxt\tE_1}\le\zz_C$. Taking maximum over all $r,q,E_1,\tE_1$ as above we obtain 
$\fa_E\le\zz_C$.

\subhead 6.8\endsubhead
Assume now that $\d_{y_*}=1$. Then $|\fI(y_*)|\ge2$. By 2.6(a) we can find $(x_*,\tx_*)\in S(y_*)$. By 2.6(c) we
have $|\fS(x_*)|\ge2$, $|\fS(\tx_*)|\ge2$. Define $(r,q)\in\NN^2$ by $x_*\in X_m^r,\tx_*\in X_m^q$. We must have 
$r+q=n$. We can find $E_1\in\cs_{W'_r}$, $\tE_1\in\cs_{W'_q}$ such that $x_*=\z'_r(E_1),\tx_*=\z'_q(\tE_1)$. As 
earlier in the proof we can assume that $E=j_{W'_r\T W'_q}^{W'_n}(E_1\bxt\tE_1)$ and we have
$$f_{E_1\bxt\tE_1}=f_{E_1}f_{\tE_1}=2^{(|\fS(x_*)|-2)/2+(|\fS(\tx_*)|-2)/2}=2^{|\fI(y_*)|-2}=\zz_C.$$

\subhead 6.9\endsubhead
Next we assume that $\fI(y_*)\ne\em$ and $\d_{y_*}=0$. By 2.6(a) we can find $(x_*,\tx_*)\in S(y_*)$. By 2.6(vii)
we have $\fS(\tx_*)=\em$. Then $|\fS(x_*)|=2|\fI(y_*)|$. Define $(r,q)\in\NN^2$ by $x_*\in X_m^r,\tx_*\in X_m^q$.
We must have $r+q=n$. We can find $E_1\in\cs_{W'_r}$, $\tE_1\in\cs_{W'_q}$ such that $x_*=\z'_r(E_1)$,
$\tx_*=\z'_q(\tE_1)$. We can assume that $E=j_{W'_r\T W'_q}^{W'_n}(E_1\bxt\tE_1)$ and we have
$$f_{E_1\bxt\tE_1}=f_{\tE_1}=2^{(|\fS(\tx_*)|-2)/2}=2^{|\fI(y_*)|-1}=\zz_C.$$
Now we assume that $\fI(y_*)=\em$. By 2.6(a) we can find $(x_*,\tx_*)\in S(y_*)$. By 2.6(b) we have
$\fS(x_*)=\em$, $\fS(\tx_*)=\em$. Define $(r,q)\in\NN^2$ by $x_*\in X_m^r,\tx_*\in X_m^q$. We must have $r+q=n$. 
We can find $E_1\in\cs_{W'_r}$, $\tE_1\in\cs_{W'_q}$ such that $x_*=\z'_r(E_1),\tx_*=\z'_q(\tE_1)$. We can assume
that $E=j_{W'_r\T W'_q}^{W'_n}(E_1\bxt\tE_1)$ and we have $f_{E_1\bxt\tE_1}=1=\zz_C$.

We see that in any case, $\fa_E=\zz_C$, proving 1.5(b2) in our case.

\subhead 6.10\endsubhead
For $g\in\Om$ let $\la g\ra$ be the subgroup of $\Om$ generated by $g$.

When $n$ is even the subgroups of $\Om$ are $\{1\},\la\o_1\ra,\la\o_2\ra,\la\o_1\o_2\ra,\Om$; when $n$ is odd the 
subgroups of $\Om$ are $\{1\},\la\o_2\ra,\Om$. 

(a) The collection of subgroups $W'_r\T S_p^{(0)}\T W'_q$ (with $r=q\ge1$) contains all subgroups in $\cp^\Om$.

(b) The collection of subgroups $W'_r\T W'_q$ contains all subgroups in $\cp^{\la\o_2\ra}$.

(c) For $n$ even, the collection in (a) together with the subgroups $W'_0\T S_p^{(\l)}\T W'_0$ (with $\l=0$ or 
$3$) contains all subgroups in $\cp^{\la\o_1\ra}$.

(d) For $n$ even, the collection in (a) together with the subgroups $W'_0\T S_p^{(\l)}\T W'_0$ (with $\l=1$ or 
$2$) contains all subgroups in $\cp^{\la\o_1\o_2\ra}$.

\subhead 6.11\endsubhead
Assume that $\tzz_C/\zz_C=4$. Then $\d_{y_*}=1$ and $|\ci|=1$ for any $\ci\in\fI(y_*)$. By 2.11 we can find $r,p$,
$x_*\in X_m^r$, $e_*\in\ce_m^p$ (with $r+p+r=n$) such that $y_*=x_*+e_*+x_*$, $(x_*,e_*+x_*)\in S(y_*)$, 
$\fS(x_*)=\fS(e_*+x_*)\ne\em$. Note that $r\ge1$. Define $E_1\in\cs_{W'_r}$ by $\z'_r(E_1)=x_*$, $E_2\in\cs_{S_p}$
by $\x_p(E_2)=e_*$. We have $E=j_{W'_r\T S_p^{(0)}\T W'_r}^{W'_n}(E_1\bxt E_2\bxt E_1)$ and
$$\align&f_{E_1\bxt E_2\bxt E_1}=2^{(|\fS(x_*)|-2)/2}2^{(|\fS(x_*)|-2)/2}=
2^{(|\fS(x_*)|-2)/2}2^{(|\fS(e_*+x_*)|-2)/2}\\&=2^{|\fI(y_*)|-2}=\zz_C.\endalign$$
We have $W'_r\T S^{0)}_p\T W'_r\in\cp^\Om$. Moreover, $E_1\bxt E_2\bxt E_1$ is $\Om$-stable. We see that $\fc_E=4$. 

\subhead 6.12\endsubhead
Conversely, assume that $\fc_E=4$. By 1.11 and 6.10(a), there exist $(r,p,q)$ as in 6.3 with $q=r\ge1$ and 
$E_1\in\cs_{W'_r}$, $E_2\in\cs_{S_p}$ such that $E=j_{W'_r\T S_p^{(0)}\T W'_r}^{W'_n}(E_1\bxt E_2\bxt E_1)$, 
$f_{E_1\bxt E_2\bxt E_1}=\zz_C$ and such that $E_1$ extends to a $W_r$-module. We set $x_*=\z'_r(E_1)\in X_m^r$, 
$e_*=\x_p(E_2)$. We have $y_*=x_*+e_*+x_*$. Since $E_1$ extends to a $W_r$-module we have $\fS(x_*)\ne\em$, hence
$\fI(y_*)\ne\em$. Thus, $\zz_C=2^{|\fI(y_*)|-1-\d_{y_*}}$, 
$2^{(|\fS(x_*)|-2+|\fS(x_*)|-2)/2}=2^{|\fI(y_*)|-1-\d_{y_*}}$ and $|\fS(x_*)|+|\fS(x_*)|=2|\fI(y_*)|+1-\d_{y_*}$. 
Since $|\fS(x_*)|+|\fS(x_*)|\le 2|\fI(y_*)|$, we have $1-\d_{y_*}\le0$ hence $\d_{y_*}=1$ and 
$|\fS(x_*)|+|\fS(x_*)|=2|\fI(y_*)|$. 

Let $E'_1=j_{S_p\T W'_r}^{W'_{p+r}}(E_2\bxt E_1)\in\cs_{W'_{p+r}}$. Then 
$E=j_{W'_r\T W'_{p+r}}^{W'_n}(E_1\bxt E'_1)$. By 1.5(b2) we have $f_{E_1\bxt E'_1}\le\zz_C$. By 1.9(b) we have 
$f_{E_2\bxt E_1}\le f_{E'_1}$. Hence $\zz_C=f_{E_1\bxt E_2\bxt E_1}\le f_{E_1\bxt E'_1}\le\zz_C$; this forces 
$f_{E_2\bxt E_1}=f_{E'_1}$. The last equality can be rewritten as
$$2^{(|\fS(x_*)|-2)/2}=2^{(|\fS(e_*+x_*)|-2)/2}$$
since $e_*+x_*=\z'_{p+r}(E'_1)$ (a consequence of 6.2(a)). Hence $|\fS(e_*+x_*)|=|\fS(x_*)|$. We have also 
$(x_*,e_*+x_*)\in S(y_*)$. Using 2.10, we see that for any $\ci\in\fI(y_*)$ we have $|\ci|=1$. Thus, 
$\tzz_C/\zz_C=4$.

Using this together with 6.11, we see that $\fc_E=4$ if and only if $\tzz_C/\zz_C=4$.

\subhead 6.13\endsubhead
Assume that $\fI(y_*)=\em$. Then $n$ is even. Define $e_*\in\NN^{m+1}$ by $y_*=x^0_*+e_*+x^0_*$. We have 
$e_*\in\ce_m^n$. Define $E_1\in\cs_{W'_0}$ by $\z'_0(E_1)=x_*^0$, $E_2\in\cs_{S_n}$ by $\x_n(E_2)=e_*$. For some 
$\l\in[0,3]$ we have $E=j_{W'_0\T S_n^{(\l)}\T W'_0}^{W'_n}(E_1\bxt E_2\bxt E_1)$, see 6.3. We have 
$f_{E_1\bxt E_2\bxt E_1}=1=\zz_C$. Note that $W'_0\T S_n^{(\l)}\T W'_0\in\cp^{\Om_1}$ where $\Om_1$ is
$\la\o_1\ra$ or $\la\o_1\o_2\ra$; moreover $E_1\bxt E_2\bxt E_1$ is $\Om_1$-stable. We see that $\fc_E\ge2$. By 
6.12 we cannot have $\fc_E=4$. Hence $\fc_E=2$.

\subhead 6.14\endsubhead
Assume that $\d_{y_*}=1$ and $|\ci|>1$ for some $\ci\in\fI(y_*)$. We have $|\fI(y_*)|\ge2$. By 2.6(a) we can find
$(x_*,\tx_*)\in S(y_*)$. By 2.6(c) we have $|\fS(x_*)|\ge2$, $|\fS(\tx_*)|\ge2$. Define $(r,q)\in\NN^2$ by 
$x_*\in X_m^r,\tx_*\in X_m^q$. We must have $r+q=n$ and $r\ge1$, $q\ge1$. We can find uniquely $E_1\in\cs_{W'_r}$,
$\tE_1\in\cs_{W'_q}$ such that $x_*=\z'_r(E_1),\tx_*=\z'_q(\tE_1)$. We have 
$E=j_{W'_r\T W'_q}^{W'_n}(E_1\bxt\tE_1)$ and
$$f_{E_1\bxt\tE_1}=f_{E_1}f_{\tE_1}=2^{(|\fS(x_*)|-2)/2+(|\fS(\tx_*)|-2)/2}=2^{|\fI(y_*)|-2}=\zz_C.$$
We have $W'_r\T W'_q\in\cp^{\la\o_2\ra}$ and $E_1\bxt\tE_1$ is $\la\o_2\ra$-stable. We see that $\fc_E\ge2$. By 
6.12 we cannot have $\fc_E=4$. Hence $\fc_E=2$.

\subhead 6.15\endsubhead
Assume that $\fc_E=2$. By 1.11 and 6.10, either (i) or (ii) below holds.

(i) there exist $(r,p,q)$ as in 6.3 with $q=r$, $\l\in[0,3]$ (with $\l=0$ unless $r=0$) and $E_1\in\cs_{W'_r}$, 
$E_2\in\cs_{S_p}$ such that $E=j_{W'_r\T S_p^{(\l)}\T W'_r}^{W'_n}(E_1\bxt E_2\bxt E_1)$, 
$f_{E_1\bxt E_2\bxt E_1}=\zz_C$;

(ii) there exist $(r,q)$ with $r+q=n$ and $E_1\in\cs_{W'_r}$, $\tE_1\in\cs_{W'_q}$ such that $E_1$ extends to a 
$W_r$-module, $\tE_1$ extends to a $W_q$-module, $E=j_{W'_r\T W'_q}^{W'_n}(E_1\bxt\tE_1)$ and
$f_{E_1\bxt\tE_1}=\zz_C$.
\nl
Assume first that (i) holds. We set $x_*=\z'_r(E_1)\in X_m^r$. If $r\ge1$ and $E_1$ extends to a $W_r$-module then
$E_1\bxt E_2\bxt E_1$ is $\Om$-stable (note that $W'_r\T S_p^{(\l)}\T W'_r\in\cp^\Om$) so that $\fc_E=4$ 
contradicting $\fc_E=2$. Thus, either $r\ge1$ and $E_1$ does not extend to a $W_r$-module or $r=0$. It follows 
that $\fS(x_*)=\em$ and $f_{E_1}=1$ so that $\zz_C=1$. Hence either $|\fI(y_*)|=0$ or $|\fI(y_*)|=2,\d_{y_*}=1$. 
In the first case we have $\tzz_C/\zz_C=2$. In the second case, using $\d_{y_*}=1$ we see that $\tzz_C/\zz_C\ge2$;
if we had $\tzz_C/\zz_C=4$ we would have $\fc_E=4$, a contradiction. Thus in both cases we have $\tzz_C/\zz_C=2$. 

Next assume that (ii) holds. We set $x_*=\z'_r(E_1)\in X_m^r$, $\tx_*=\z'_q(\tE_1)\in\tX_m^q$. We have 
$y_*=x_*+\tx_*$. Since $E_1$ extends to a $W_r$-module and $\tE_1$ extends to a $W_q$-module we have 
$|\fS(x_*)|\ge2$, $|\fS(\tx_*)|\ge2$. Hence
$$2^{|\fI(y_*)|-2}\le\zz_C=f_{E_1\bxt\tE_1}=f_{E_1}f_{\tE_1}=2^{(|\fS(x_*)|-2)/2+(|\fS(\tx_*)|-2)/2}\le 
2^{|\fI(y_*)|-2}.$$
It follows that $2^{|\fI(y_*)|-2}=\zz_C$ so that $\d_{y_*}=1$. This implies that $\tzz_C/\zz_C\ge2$; if we had 
$\tzz_C/\zz_C=4$ we would have $\fc_E=4$, a contradiction. Thus we have $\tzz_C/\zz_C=2$. 

Using this together with 6.13, 6.14, we see that $\fc_E=2$ if and only if $\tzz_C/\zz_C=2$.

\subhead 6.16\endsubhead
By 6.12, we have $\fc_E=4$ if and only if $\tzz_C/\zz_C=4$. By 6.15, we have $\fc_E=2$ if and only if 
$\tzz_C/\zz_C=2$. Since $\fc_E\in\{1,2,4\}$ and $\tzz_C/\zz_C\in\{1,2,4\}$ we see that $\fc_E=\tzz_C/\zz_C$; this 
proves 1.5(b3) in our case.

\head 7. Exceptional types\endhead
\subhead 7.1\endsubhead
In this section we assume that $G$ is an exceptional group. For each type we give a table with rows indexed by the
unipotent conjugacy classes in $G$ in which the row corresponding to $C\in\cx$ has four entries:
$$\r_C\qquad\bb_C\qquad a\T a'\quad(J,E_1)$$
where $a=\zz_C$, $a'=\tzz_C/\zz_C$ and $(J,E_1)$ is an example of an element of $\cz_E$ ($E=\r_C$) such that 
$f_{E_1}=\zz_C$ and $|\Om_{J,E_1}|=\tzz_C/\zz_C$. (When $\Om=\{1\}$ we have $a'=1$ and we write $a$ instead of 
$a\T a'$). 
We specify an irreducible representation $E_1$ of a Weyl group either by using the notation of \cite{\OR, Ch.4} 
(for type $E_6,E_7,E_8$) or by specifying its degree. The representation is then determined by its $b_{E_1}$ which
equals $\bb_C$ in the table or (in the case of $G_2$, $F_4$) by other information in the same row of the table. On
the other hand, $\e$ always denotes the sign representation. In a pair $(J,E_1)$, $J$ is any subset of $\tI$ such 
that $\WW_J$ has the specified type; in addition, for type $F_4$, we denote by $A_2$ (resp. $A'_2$) a subset $J$ 
of $\tI$ such that $\WW_J$ is of type $A_2$ and is contained (resp. not contained) in a parahoric subgroup of type
$B_4$).

The group $\Om$ is $\{1\}$ for types $G_2,F_4$ and is a cyclic group of order
$9-n$ for type $E_n (n=6,7,8)$.

Type $G_2$     
$$\allowdisplaybreaks\alignat4
\r_C&&\qquad\bb_C&&\qquad a\T a'&&\quad(J,E_1)\\
1&&\qquad0 &&\qquad   1 && \qquad       (\em,1)           \\
2&&\qquad1 && \qquad  6 && \qquad       (G_2,2)          \\
2&&\qquad2 &&\qquad  1  &&\qquad        (A_1A_1,\e)       \\
1&&\qquad3 &&\qquad  1  &&\qquad        (A_2,\e) \\        
1&&\qquad6 &&\qquad  1  &&\qquad         (G_2,\e)            
\endalignat$$           

Type $F_4$
$$\allowdisplaybreaks\alignat4
\r_C&&\qquad\bb_C&&\qquad a\T a'&&\quad(J,E_1)\\
1&&\qquad0 &&\qquad    1&&\qquad   (\em,1)   \\
4&&\qquad1 &&\qquad    2&&\qquad       (F_4,4) \\     
9&&\qquad2 &&\qquad    2&&\qquad       (F_4,9) \\      
8&&\qquad3 &&\qquad    1&&\qquad      (A_2,\e)  \\      
8&&\qquad3 &&\qquad    1 &&\qquad    (A'_2,\e) \\         
12&&\qquad4 &&\qquad   24 &&\qquad    (F_4,12)\\       
16&&\qquad5 &&\qquad  2 &&\qquad   (C_3A_1,3\bxt\e)\\  
9&&\qquad6 &&\qquad    2&&\qquad    (B_4,6)\\       
6&&\qquad6 &&\qquad    1&&\qquad       (A_2A'_2,\e)\\    
4&&\qquad7 &&\qquad    1&&\qquad       (A_3A_1,\e)\\    
8&&\qquad9 &&\qquad    1&&\qquad       (C_3,\e)\\        
8&&\qquad9 &&\qquad    2&&\qquad        (B_4,4)\\  
9&&\qquad10&&\qquad    1&&\qquad     (C_3A_1,\e) \\     
4&&\qquad13&&\qquad    2&&\qquad      (F_4,4)\\         
2&&\qquad16&&\qquad    1&&\qquad       (B_4,\e)\\       
1&&\qquad24&&\qquad    1&&\qquad       (F_4,\e)
\endalignat$$           
              
Type $E_6$
$$\allowdisplaybreaks\alignat4
\r_C&&\qquad\bb_C&&\qquad a\T a'&&\quad(J,E_1)\\
1_p&&\qquad0 &&\qquad   1\T3  &&\qquad        (\em,1)           \\
6_p&&\qquad1 &&\qquad   1\T3  &&\qquad      (D_4,4) \\          
20_p&&\qquad2&&\qquad   1\T1&&\qquad       (E_6,20_p)\\       
30_p&&\qquad3&&\qquad   2\T3 &&\qquad     (D_4,8)\\    
15_q&&\qquad4 &&\qquad  1\T3 &&\qquad     (A_1A_1A_1A_1,\e)\\ 
64_p&&\qquad4 &&\qquad  1\T1&&\qquad         (E_6,64_p)       \\
60_p&&\qquad5&&\qquad   1\T1&&\qquad         (E_6,60_p)    \\
24_p&&\qquad6&&\qquad   1\T1&&\qquad          (E_6,24_p)       \\
81_p&&\qquad6 &&\qquad  1\T1&&\qquad         (E_6,81_p)         \\
80_s&&\qquad7&&\qquad   6\T1&&\qquad       (E_6,80_s)           \\
60_s&&\qquad8&&\qquad   1\T1&&\qquad          (A_3A_1A_1,\e)    \\
10_s&&\qquad9&&\qquad   1\T3&&\qquad       (A_2A_2A_2,\e)    \\
81'_p&&\qquad10&&\qquad  1\T1&&\qquad          (E_6,81'_p)          \\
60'_p&&\qquad11&&\qquad  1\T1&&\qquad         (E_6,60'_p)      \\
24'_p&&\qquad12&&\qquad  1\T3&&\qquad       (D_4,\e)          \\
64'_p&&\qquad13&&\qquad  1\T1&&\qquad         (E_6,64'_p)      \\
30'_p&&\qquad15&&\qquad  2\T1&&\qquad       (E_6,30'_p)      \\
15'_q&&\qquad16&&\qquad  1\T1&&\qquad          (A_5A_1,\e)\\     
20'_p&&\qquad20&&\qquad  1\T1&&\qquad     (E_6,20'_p)\\      
6'_p&&\qquad25 &&\qquad 1\T1&&\qquad       (E_6,6'_p)            \\
1'_p&&\qquad36  &&\qquad1\T1&&\qquad        (E_6,\e)           
\endalignat$$

Type $E_7$
$$\allowdisplaybreaks\alignat4
\r_C&&\qquad\bb_C&&\qquad a\T a'&&\quad(J,E_1)\\
1_a&&\qquad0 &&\qquad   1\T 2 &&\qquad     (\em,1)     \\
7'_a&&\qquad1&&\qquad    1\T 2 &&\qquad    (E_6,6_p)            \\
27_a&&\qquad2 &&\qquad  1\T 2 &&\qquad     (E_6,20_p)        \\
56'_a&&\qquad3 &&\qquad  2\T2 &&\qquad      (E_6,30_p)     \\
21'_b&&\qquad3 &&\qquad  1\T1  &&\qquad      (E_7,21'_b)    \\
120_a&&\qquad4 &&\qquad 2\T1 &&\qquad    (E_7,120_a)  \\
35_b&&\qquad4  &&\qquad 1\T2 &&\qquad    (A_7,14) \\   
189'_b&&\qquad5 &&\qquad2\T2 &&\qquad     (A_1D_4A_1,\e\bxt8\bxt\e) \\
105_b&&\qquad6 &&\qquad 1\T1 &&\qquad        (E_7,105_b)\\ 
210_a&&\qquad6 &&\qquad  1\T2 &&\qquad      (A_7,35)\\          
168_a&&\qquad6 &&\qquad  1\T2 &&\qquad    (A_7,56)\\            
315'_a&&\qquad7&&\qquad   6\T2&&\qquad     (E_6,80_s)\\        
189'_c&&\qquad7&&\qquad   1\T1&&\qquad        (E_7,189'_c)\\    
405_a&&\qquad8&&\qquad   2\T1&&\qquad      (E_7,405_a)\\  
280_b&&\qquad8&&\qquad   1\T2&&\qquad      (A_7,56)\\           
70'_a&&\qquad9&&\qquad   1\T2&&\qquad      (A_2A_2A_2,\e)\\   
216'_a&&\qquad9&&\qquad   1\T1&&\qquad        (D_6A_1,30\bxt\e)      \\
378'_a&&\qquad9&&\qquad  1\T2 &&\qquad    (A_7,70)\\          
420_a&&\qquad10&&\qquad  2\T1 &&\qquad     (E_7,420_a)\\    
210_b&&\qquad10&&\qquad  1\T1&&\qquad         (E_7,210_b)\\          
512'_a&&\qquad11&&\qquad  2\T1&&\qquad     (E_7,512'_a)\\                
105_c&&\qquad12&&\qquad  1\T2&&\qquad      (D_4,\e)             \\
84_a&&\qquad12 &&\qquad 1\T2 &&\qquad       (A_7,14)\\                        
420'_a&&\qquad13 &&\qquad  2\T1&&\qquad       (D_6,24)      \\         
210_b&&\qquad13  &&\qquad1\T2 &&\qquad       (A_3A_3A_1,\e)            \\
378'_a&&\qquad14 &&\qquad 2\T1&&\qquad  (D_6A_1,24\bxt\e)  \\       
105'_c&&\qquad15 &&\qquad 1\T1  &&\qquad   (A_5A_2,\e\bxt1)\\             
405'_a&&\qquad15 &&\qquad 2\T2 &&\qquad  (E_6,30'_p)\\           
216_a&&\qquad16  &&\qquad1\T2 &&\qquad    (A_7,20)\\      
315_a&&\qquad16  &&\qquad6\T1 &&\qquad        (E_7,315_a)\\                  
280'_b&&\qquad17  &&\qquad1\T1 &&\qquad         (D_6A_1,15\bxt\e) \\
70_a&&\qquad18  &&\qquad1\T1 &&\qquad              (A_5A_2,\e)             \\
189_c&&\qquad20 &&\qquad1\T2 &&\qquad        (E_6,20'_p)\\               
210'_a&&\qquad21 &&\qquad 1\T1&&\qquad        (E_7,210'_a)               \\ 
168'_a&&\qquad(21 &&\qquad    1\T1&&\qquad        (E_7,168'_a)            \\
105'_b&&\qquad21 &&\qquad1\T2 &&\qquad         (A_7,7)         \\
189_b&&\qquad22   &&\qquad1\T1&&\qquad         (E_7,189_b)    \\
120'_a&&\qquad25  &&\qquad   2\T1&&\qquad     (E_7,120'_a)                  \\
15_a&&\qquad28  &&\qquad   1\T2 &&\qquad        (A_7,\e)                   \\
56_a&&\qquad30  &&\qquad   2\T1 &&\qquad       (E_7,56_a)                      \\
35'_b&&\qquad31 &&\qquad   1\T1 &&\qquad        (D_6A_1,\e)             \\
21_b&&\qquad36 &&\qquad    1\T2&&\qquad       (E_6,\e)                \\
27'_a&&\qquad37 &&\qquad   1\T1 &&\qquad        (E_7,27'_a)                    \\
7_a&&\qquad46  &&\qquad   1\T1&&\qquad           (E_7,7_a)               \\
1'_a&&\qquad63 &&\qquad    1\T1  &&\qquad        (E_7,\e)                 
\endalignat$$

Type $E_8$
$$\allowdisplaybreaks\alignat4
\r_C&&\qquad\bb_C&&\qquad a\T a'&&\quad(J,E_1)\\
1_x&&\qquad0 &&\qquad   1 &&\qquad                (\em,1)           \\
8_z&&\qquad1 &&\qquad     1&&\qquad              (E_8,8_z)\\          
35_x&&\qquad2&&\qquad     1&&\qquad              (E_8,35_x)\\       
112_z&&\qquad3&&\qquad   2 &&\qquad             (E_8,112_z) \\
84_x&&\qquad4&&\qquad 1  &&\qquad             (E_7A_1,21'_b\bxt\e)    \\
210_x&&\qquad4&&\qquad   2&&\qquad              (E_8,210_x)\\          
560_z&&\qquad5&&\qquad   2 &&\qquad             (E_7A_1,120_a\bxt\e)\\
567_x&&\qquad6&&\qquad   1 &&\qquad             (E_8,567_x)\\           
700_x&&\qquad6&&\qquad   2 &&\qquad             (E_8,700_x)\\         
400_x&&\qquad7 &&\qquad1   &&\qquad           (A_2A_1A_1A_1A_1,\e)\\ 
1400_z&&\qquad7 &&\qquad 6 &&\qquad              (E_6,80_s)            \\
1400_x&&\qquad8 &&\qquad 6 &&\qquad              (E_8,1400_x)  \\  
1344_x&&\qquad8 &&\qquad 1 &&\qquad               (E_7A_1,189'_c\bxt\e)  \\
448_z&&\qquad9 &&\qquad 1  &&\qquad                (A_2A_2A_2,\e)      \\
3240_z&&\qquad9 &&\qquad 2  &&\qquad                (E_7A_1,405_a\bxt\e)   \\
2240_x&&\qquad10 &&\qquad6  &&\qquad       (E_6A_2,80_s\bxt\e)\\    
2268_x&&\qquad10 &&\qquad   2 &&\qquad (E_8,2268_x)         \\
4096_x&&\qquad11 &&\qquad 2 &&\qquad (E_7,512'_a)\\       
1400_z&&\qquad11 &&\qquad 1 &&\qquad  (E_7A_1,210_b\bxt\e)\\   
525_x&&\qquad12 &&\qquad  1 &&\qquad       (D_4,\e)\\              
4200_x&&\qquad12&&\qquad    2&&\qquad    (E_8,4200_x) \\
 972_x&&\qquad12&&\qquad  1 &&\qquad      (A_3A_3,\e) \\         
2800_z&&\qquad13&&\qquad    2&&\qquad      (E_8,2800_z)\\               
4536_z&&\qquad13&&\qquad   2 &&\qquad    (D_8,560)\\    
6075_x&&\qquad14&&\qquad    2 &&\qquad    (D_8,280) \\
2835_x&&\qquad14&&\qquad   1  &&\qquad     (A_4A_2A_1,\e)\\         
4200_z&&\qquad15&&\qquad   1  &&\qquad   (A_5,\e)\\             
5600_z&&\qquad15 &&\qquad 2  &&\qquad       (E_6,30'_p)\\           
4480_y&&\qquad16 &&\qquad120 &&\qquad     (E_8,4480_y)               \\
3200_x&&\qquad16 &&\qquad  1 &&\qquad    (A_5A_1,\e)         \\
7168_w&&\qquad17 &&\qquad  6&&\qquad  (E_7A_1,315_a\bxt\e)\\
4200_y&&\qquad18 &&\qquad   2 &&\qquad       (D_8,252)              \\
3150_y&&\qquad18 &&\qquad    2 &&\qquad     (E_6A_2,30'_p\bxt\e)       \\
2016_w&&\qquad19 &&\qquad    1 &&\qquad    (A_5A_2A_1,\e)       \\
1344_w&&\qquad19 &&\qquad  1  &&\qquad       (D_5A_3,5\bxt\e)         \\
2100_y&&\qquad20 &&\qquad     1 &&\qquad     (D_5,\e)              \\
420_y&&\qquad20  &&\qquad     1 &&\qquad      (A_4A_4,\e)           \\
5600'_z&&\qquad21 &&\qquad2 &&\qquad         (E_8,5600'_z)             \\
4200'_z&&\qquad21 &&\qquad   2&&\qquad  (D_8,224)         \\
3200'_x&&\qquad22 &&\qquad   1 &&\qquad      (E_7A_1,168_a\bxt\e) \\   
6075'_x&&\qquad22 &&\qquad1 &&\qquad        (E_8,6075'_x)            \\
2835'_x&&\qquad22  &&\qquad   1 &&\qquad     (A_6A_1,\e)          \\
4536'_z&&\qquad23  &&\qquad   1 &&\qquad     (D_5A_2,\e)         \\
4200'_x&&\qquad24  &&\qquad 2 &&\qquad     (E_8,4200'_x)         \\
2800'_z&&\qquad25 &&\qquad  2  &&\qquad      (E_7,120'_a)              \\
4096'_x&&\qquad26 &&\qquad    2 &&\qquad     (E_8,4096'_x)          \\
840'_x&&\qquad26  &&\qquad  1  &&\qquad      (D_5A_3,\e)          \\
700'_x&&\qquad28  &&\qquad 1  &&\qquad       (A_7,\e)                   \\
2240'_x&&\qquad28 &&\qquad  2 &&\qquad     (E_8,2240'_x)                \\
1400'_z&&\qquad29 &&\qquad  1 &&\qquad     (A_7A_1,\e)                  \\
2268'_x&&\qquad30 &&\qquad     2&&\qquad         (E_7,56_a)                \\
3240'_z&&\qquad31 &&\qquad    2&&\qquad     (E_7A_1,56_a\bxt\e)      \\
1400'_x&&\qquad32 &&\qquad    6&&\qquad    (E_8,1400'_x)            \\
1050'_x&&\qquad34 &&\qquad     1 &&\qquad      (D_8,28)                \\
 525'_x&&\qquad36 &&\qquad   1  &&\qquad        (E_6,\e))                \\
175'_x&&\qquad36  &&\qquad      1&&\qquad       (A_8,\e)                 \\
1400'_z&&\qquad37 &&\qquad  6 &&\qquad           (E_8,1400'_z)                \\
1344'_x&&\qquad38 &&\qquad  1 &&\qquad   (E_7A_1,27'_a\bxt\e)              \\
448'_z&&\qquad39  &&\qquad   1 &&\qquad          (E_6A_2,\e)\\               
700'_x&&\qquad42  &&\qquad    2 &&\qquad         (E_8,700'_x)\\                       
400'_z&&\qquad43  &&\qquad    1 &&\qquad   (D_8,8)               \\
567'_x&&\qquad46  &&\qquad  1  &&\qquad           (E_7,7_a)                       \\
560'_z&&\qquad47  &&\qquad   1 &&\qquad      (E_7A_1,7_a\bxt\e)          \\
210'_x&&\qquad52  &&\qquad     2&&\qquad         (E_8,210'_x)                   \\
50'_x&&\qquad56   &&\qquad     1&&\qquad      (D_8,\e)                        \\
112'_z&&\qquad63  &&\qquad    2&&\qquad   (E_8,112'_z)               \\
84'_x&&\qquad64  &&\qquad    1 &&\qquad       (E_7A_1,\e)                    \\
 35'_x&&\qquad74  &&\qquad    1 &&\qquad        (E_8,35'_x)                  \\
8'_z&&\qquad91   &&\qquad    1  &&\qquad      (E_8,8'_z)\\                          
 1'_x&&\qquad120 &&\qquad    1  &&\qquad       (E_8,\e)                  
\endalignat$$

\head Index \endhead
0.1: $\cx,\r_C,\bb_C,\zz_C,\tzz_C,\tcs_\WW$

1.1: $\Irr(W),f_E,a_E,b_E,\cs_W,\Irr(W)^\da$

1.2: $\car,\tI,\tca,\WW,\WW_J,s_i,\Om$

1.3: $j_{\WW_J}^\WW(E_1),\bcs_\WW,\cz_E,\fa_E,\cz_E^\sp,\Om_{J,E_1},\fc_E$

1.4: $G$

1.9: $\cp'$

1.11: $\cp^{\tiO}$

2.9: $\ce_m$

4.1: $W_n$

4.2: $U_k^n,\z_n$

4.6: $\t$

5.1: $W'_n$

5.2: $V_k^n,{}^\da V_k^n,\z'_n,\tiz_n$

5.4: $\ti\t$

6.4: $\t$


\widestnumber\key{AL}
\Refs
\ref\key\AL\by D.Alvis\paper Induce/restrict matrices for exceptional Weyl groups\jour math.RT/0506377\endref
\ref\key\ALL\by D.Alvis and G.Lusztig\paper On Springer's correspondence for simple groups of type $E_n$
($n=6,7,8$)\jour Math. Proc. Camb. Phil. Soc.\vol92\yr1982\pages65-78\endref
\ref\key\CL\by G.Lusztig\paper Irreducible representations of finite classical groups\jour Inv.Math.\vol43\yr1977
\pages125-175\endref
\ref\key\SPE\by G.Lusztig\paper A class of irreducible representations of a Weyl group\jour Proc. Kon. Nederl.
Akad. (A)\vol82\yr1979\pages323-335\endref
\ref\key\UC\by G.Lusztig\paper Unipotent characters of the symplectic and odd
orthogonal groups over a finite field\jour Invent.math.\vol64\yr1981\pages263-296\endref
\ref\key\OR\by G.Lusztig\book Characters of reductive groups over a finite field\bookinfo Ann.Math.Studies\vol107
\publ Princeton U.Press\yr1984\endref 
\ref\key\ICC\by G.Lusztig\paper Intersection cohomology complexes on a reductive group\jour Invent.Math.\vol75
\yr1984\pages205-272\endref 
\ref\key\USU\by G.Lusztig\paper A unipotent support for irreducible representations\jour Adv.in Math.\vol 94
\yr1992\pages139-179\endref
\ref\key\SH\by T.Shoji\paper On the Springer representations of Weyl groups of classical algebraic groups\jour 
Comm.in Alg.\vol7\yr1979\pages1713-1745,2027-2033\endref
\ref\key\SHO\by T.Shoji\paper On the Springer representations of Chevalley groups of type $F_4$\jour Comm.in Alg.
\vol8\yr1980\pages409-440\endref
\ref\key\SP\by T.A.Springer\paper Trigonometric sums, Green functions of finite groups and representations of Weyl
groups\jour Invent.Math.\vol36\yr1976\pages173-207\endref
\endRefs
\enddocument